\documentclass[11pt]{article}

\usepackage{graphicx}
\usepackage{amsmath}
\usepackage{amsfonts}
\usepackage{color}

\newcommand{\be}{\begin{equation}}
\newcommand{\ee}{\end{equation}}
\newcommand{\beqn}{\begin{eqnarray}}
\newcommand{\eeqn}{\end{eqnarray}}
\newcommand{\beqns}{\begin{eqnarray*}}
\newcommand{\eeqns}{\end{eqnarray*}}
\newcommand{\lkr}{\left(}
\newcommand{\lkv}{\left[}
\newcommand{\rkv}{\right]}
\newcommand{\rkr}{\right)}
\newcommand{\lfi}{\left\{}
\newcommand{\rfi}{\right\}}

\newcommand{\fr}[1]{(\ref{#1})}

\newcommand{\del}{\delta}

\newcommand{\eps}{\varepsilon}

\newcommand{\om}{\omega}
\newcommand{\lam}{\lambda}

\newcommand{\sig}{\sigma}

\newcommand{\Lam}{\Lambda}

\newcommand{\EE}{\ensuremath{{\mathbb E}}}

\newcommand{\II}{\ensuremath{{\mathbb I}}}

\newcommand{\PP}{\ensuremath{{\mathbb P}}}

\newcommand{\RR}{\ensuremath{{\mathbb R}}}
\newcommand{\WW}{\mbox{\mathversion{bold}${\ensuremath{\mathbb W}}$}}

\newcommand{\card}{\mbox{card}}

\newcommand{\Tr}{{\rm Tr}}
\newcommand{\Vect}{{\rm Vec}}
\newcommand{\block}{\mbox{block}}

\newtheorem{theorem}{Theorem}
\newtheorem{lemma}{Lemma}
\newtheorem{corollary}{Corollary}
\newtheorem{proposition}{Proposition}
\newtheorem{remark}{Remark}
\newtheorem{example}{Example}

\newcommand{\bA}{\mathbf{A}}
\newcommand{\boB}{\mathbf{B}}

\newcommand{\bI}{\mathbf{I}}
\newcommand{\bQ}{\mathbf{Q}}

\newcommand {\bW}{\mathbf{W}}
\newcommand {\bY}{\mathbf{Y}}
 \newcommand {\bZ}{\mathbf{Z}}
 
\newcommand{\ba}{\mathbf{a}}
\newcommand{\bbb}{\mathbf{b}}
\newcommand{\bc}{\mathbf{c}}

\newcommand{\bt}{\mathbf{t}}

\newcommand{\bv}{\mathbf{v}}
\newcommand{\bq}{\mathbf{q}}
\newcommand{\br}{\mathbf{r}}

\newcommand {\bone}{\mathbf{e}}
\newcommand{\bof}{\mathbf{f}}
\newcommand{\bzero}{\mathbf{0}}

\newcommand{\calL}{{\cal L}}
\newcommand{\calJ}{{\cal J}}

\newcommand {\brho}{\mbox{\mathversion{bold}$\rho$}}
\newcommand {\bphi}{\mbox{\mathversion{bold}$\phi$}}
\newcommand {\balpha}{\mbox{\mathversion{bold}$\alpha$}}
\newcommand {\bxi}{\mbox{\mathversion{bold}$\xi$}}
\newcommand {\bnu}{\mbox{\mathversion{bold}$\nu$}}

\newcommand {\bom}{\mbox{\mathversion{bold}$\om$}}
\newcommand {\bzeta}{\mbox{\mathversion{bold}$\zeta$}}

\newcommand {\bOmega}{\mbox{\mathversion{bold}$\Omega$}}
\newcommand {\bPhi}{\mbox{\mathversion{bold}$\Phi$}}

\newcommand {\bSigma}{\mbox{\mathversion{bold}$\Sigma$}}
\newcommand {\bSig}{\mbox{\mathversion{bold}$\Sigma$}}

\newcommand {\cOm}{\Theta}

\newcommand{\Sum}{\sum_{i=1}^n}

\newcommand{\hba}{\hat{\ba}}
\newcommand{\hbA}{\widehat{\bA}}
\newcommand {\hbSigma}{\mbox{\mathversion{bold}$\widehat{\bSigma}$}}
\newcommand {\hbSig}{\mbox{\mathversion{bold}$\widehat{\bSigma}$}}

\newcommand{\hbf}{\hat{\bof}}
\newcommand{\tbf}{\tilde{\bof}}
\newcommand{\tbof}{\tilde{\bof}}
\newcommand{\hf}{\hat{f}}
\newcommand{\ha}{\hat{a}}

\newcommand {\hdel}{\hat{\del}}

\newcommand {\tbom}{\tilde{\bom}}
\newcommand {\tom}{\tilde{\om}}

\newcommand{\Umu}{U_{\mu}}
\newcommand{\Cmu}{C_{\mu}}
\newcommand{\wpmu}{{\cal W}_{\mu}}
\newcommand{\wpmun}{\wpmu^{\otimes n}}
\newcommand{\calF}{{\cal F}}
\newcommand{\tcalF}{\tilde{\calF}}

\newcommand{\nuj}{{\nu_j}}
\newcommand{\rj}{{r_j}}
\newcommand{\rmax}{{r_{max}}}

\newcommand{\omin}{\omega_{\min}}
\newcommand{\omax}{\omega_{\max}}
\newcommand{\omins}{\omega^*_{\min}}
\newcommand{\omaxs}{\omega^*_{\max}}

\newcommand{\phimin}{\phi_{\min}}
\newcommand{\phimax}{\phi_{\max}}
\newcommand{\lamin}{\lam_{\min}}
\newcommand{\lamax}{\lam_{\max}}

\def\keywords{\vspace{.5em}
{\textit{Keywords}:\,\relax%
}}
\def\endkeywords{\par}

\def\subjclass{\vspace{.5em}
{\textit{AMS 2000 subject classification}:\,\relax%
}}
\def\endsubjclass{\par}

\begin{document}

\title{\bf { Sparse high-dimensional varying coefficient model: non-asymptotic minimax study }}
   \author{{\em Olga Klopp},   CREST and MODAL'X, University Paris Ouest, 92001 Nanterre, France\\
{\em Marianna Pensky}, University of Central Florida, Orlando, FL 32816, USA}
  
\date{}

 \maketitle
 
\begin{abstract}
 The objective of the present paper is to develop a minimax theory for the varying coefficient model in a non-asymptotic setting.
 We  consider a  high-dimensional sparse varying coefficient model where  
only few of the covariates are present and only some of those  covariates are time dependent.  
Our analysis allows the time dependent covariates to have different degrees of smoothness and  to be spatially inhomogeneous.  
We develop the minimax lower bounds for the quadratic risk and construct an adaptive estimator which attains those 
lower bounds within   a constant (if all time-dependent covariates are spatially homogeneous) or logarithmic factor of the number of observations. 
\end{abstract}
%

\keywords {varying coefficient model, sparse model, minimax optimality}
\endkeywords

\subjclass{62H12, 62J05, 62C20 }
\endsubjclass
\section{Introduction}
\label{sec:intro}
\setcounter{equation}{0}

One of the fundamental tasks in statistics is to characterize the relationship between a set of covariates and a response variable. 
In the present paper we study the varying coefficient model which is  commonly used for describing time-varying covariate effects. 
It  provides a more  flexible approach than the classical linear regression model and is often used to 
analyze the data measured repeatedly over time.

Since its introduction by Cleveland, Grosse and Shyu \cite{cleveland_VCM} and 
Hastie and Tibshirani \cite {hastie_VCM} many methods for estimation and inference  in the varying coefficient model  have been developed
(see, e.g.,  \cite{wu,hoover,fan_1999,kauermann} for the kernel-local polynomial smoothing approach, \cite{huang_2002,huang_2004,huang_shen} 
for the polynomial spline approach, \cite{hastie_VCM,hoover,chiang} for the smoothing spline approach and \cite{fan_2008} 
for a detailed discussion of the existing methods and possible applications). In the last five years, the varying coefficient 
model received a great deal of attention.  For example, 
Wang {\it et al.}  \cite{wang} proposed a new procedure based on a local rank estimator; 
Kai {\it et al.}  \cite{kai} introduced a semi-parametric quantile regression procedure and studied an 
effective variable selection procedure. Lee {\it et al.}  \cite{lee} extended the model to the case
when the response is related to the covariates via a link function while Zhu {\it et al.} \cite{zhu} studied the 
multivariate version of the model.   Existing methods  typically provide asymptotic evaluation of the precision of the estimation procedure 
under the assumption that the number of observations tends to infinity and is larger than the dimension of the problem.

  Recently few authors consider still asymptotic but high-dimensional approach to the problem.
Wei  {\it et al.} \cite{wei}  applied group Lasso for variable selection,  
while Lian \cite{liansinica} used extended Bayesian information criterion.
 Fan   {\it et al.} \cite{fan_2013} applied nonparametric independence screening.  
Their results were extended by Lian and Ma \cite{lianma}  to include rank selection in addition to variable selection.

 One important aspect  that has  not been well studied in the existing literature is the non-asymptotic approach to the  estimation, 
prediction and variables selection in the varying coefficient model. Here, we refer to the situation where
both the number of unknown parameters  and the number of observations are large and the former may be of much higher dimension than latter. 
 The only reference that we are aware of in this setting,  is the recent paper by Klopp and Pensky \cite{klopp}. Their method is based on some recent developments in the matrix estimation problem.  Some interesting questions arise in this non-asymptotic setting.  One of them is the fundamental question of 
the minimax optimal rates of convergence. The minimax risk characterizes the essential statistical difficulty 
of the problem. It also captures the interplay between different parameters in the model. To the best of our knowledge, 
our paper presents the first \textit{non-asymptotic minimax study} 
of the sparse heterogeneous varying coefficient model.

Modern technologies  produce   very high dimensional data sets and, hence, stimulate an enormous interest in variable 
selection and estimation under a sparse scenario.  In such scenarios, penalization-based methods are particularly attractive.
Significant progress has been made in understanding the statistical properties of these methods.  
For example, many authors have studied the variable selection, 
estimation and prediction properties of the LASSO in high-dimensional setting  (see, e.g., 
 \cite{bickel_ritov_tsybakov}, \cite{bunea_tsybakov_1}, \cite{bunea_tsybakov_2}, \cite{van-geer}).  
A related LASSO-type procedure is the group-LASSO, where the covariates are assumed to 
be clustered in groups  (see, for example, \cite{yuan,bach_group,chesneau,meier_1,meier_2,lounici_pontil_tsybakov}, 
 and references therein).

In the present paper, we also consider the case when the solution is sparse, in particular, 
only few of the covariates are present and only some of them 
are time dependent.   This setup is close to  the one studied in a recent paper of Liang \cite{liansinica}. 
One important difference, however,  is that  in \cite{liansinica}, the estimator is not adaptive to 
the smoothness of the time dependent covariates. 
In addition, Liang \cite{liansinica}  assumes that all time dependent 
covariates have the same degree of smoothness and are spatially homogeneous. On the contrary, we consider a much more flexible and realistic 
scenario where the time dependent covariates possibly have different degrees of smoothness and may be spatially inhomogeneous.

In order to construct a minimax optimal  estimator, we introduce the block Lasso which can be viewed as  a version of the  group LASSO. However, 
unlike in group LASSO, where the groups occur naturally, the blocks in block LASSO are driven by the need to reduce the variance 
as it is done, for example,  in block thresholding.
Note that our estimator does not require the knowledge  which of the covariates are 
indeed present and which are time dependent. It adapts to sparsity, to heterogeneity of 
the time dependent covariates and to their possibly spatial inhomogeneous nature. In order to ensure the optimality,
we derive minimax lower bounds for the risk and show that our estimator attains those bounds within a constant (if all time-dependent 
covariates are spatially homogeneous) or logarithmic factor of the number of observations.
 The analysis is carried out under the   flexible  assumption that the noise variables are sub-Gaussian. In addition,  it does not 
require that the elements of the dictionary are uniformly bounded. 

The rest of the paper is organized  as follows. Section~\ref{sec:formulation} provides 
formulation of the problem while Section~\ref{sec:tensor} lays down a tensor approach to estimation.
Section~\ref{sec:assump} introduces  notations and assumptions on the model and  provides a discussion of the assumptions.  
Section~\ref{sec:est-err-bou} describes the block thresholding LASSO estimator, evaluates the non-asymptotic lower and upper bounds 
for the risk,   both in probability and in the mean squared risk sense, and ensures optimality of the constructed estimator.  
Section~\ref{sec:exam_discussion} presents examples of estimation when assumptions of the apper are satisfied. 
Section~\ref{sec:proofs} contains proofs of the statements formulated in the paper.


\subsection{Formulation of the problem}
\label{sec:formulation}

Let $(\bW_i,t_i,Y_i)$, $i=1,\dots,n$ be sampled independently from the   varying coefficient model 
\begin{equation}\label{vcm}
Y=\bW^{T}\bof(t)+\sigma \xi.
\end{equation} 
 Here, $\bW \in  \RR^{p}$ are random vectors of predictors,
$\bof(\cdot)=\left (f_1(\cdot),\dots,f_p(\cdot)\right )^{T}$ is an unknown vector-valued  
function of regression coefficients and $t\in[0,1]$ is a random variable    
with the unknown density function $g$.   We assume that  $\bW$ and $t$ are independent. 
 The   noise variable  $\xi$ is independent of $W$ and $t$ 
 and is such that 
 $\EE(\xi)=0$ and $\EE(\xi^{2})=1$, $\sigma>0$ denotes the known  noise level.

The goal is to estimate   vector function $f(\cdot)$ on the basis of observations  $(\bW_i,t_i,Y_i)$, $i=1,\dots,n$.

In order to estimate $\bof$, following Klopp and Pensky (2013), we expand it over a   basis 
$(\phi_l(\cdot)), {l=0,1, \dots,\infty}$,   in $L_2  ([0,1])$   with $\phi_0(t) =1$.
%
Expansion   of the functions $f_j(\cdot)$ over the basis, for any $t \in [0,1]$, yields
\be \label{func_expan}
f_j(t) = \sum_{l=0}^L a_{jl} \phi_l(t) + \rho_j(t) \quad \mbox{with} 
\quad \rho_j(t) = \sum_{l=L+1}^\infty  a_{jl} \phi_l(t).
\ee
If  $\bphi(\cdot)=(\phi_0(\cdot),\dots, \phi_L(\cdot))^{T}$  and  
$\bA_{0}$ denotes a matrix of coefficients  with elements $\bA_{0}^{(l,j)}= a_{jl}$, 
then relation \fr{func_expan} can be re-written as $\bof(t) = \bA_{0}^T \bphi(t) + \brho(t)$, where 
 $\brho(t) = (\rho_1(t), \cdots, \rho_p(t))^T$. 
Combining formulae \fr{vcm} and \fr{func_expan}, we obtain the following model for 
observations $(\bW_i,t_i,Y_i)$, $i=1,\dots,n$:
\be \label{obs_mod}
Y_i = \Tr(\bA_{0}^T \bphi(t_i) \bW_i^T) + \bW_i^T \brho(t_i) + \sig \xi_i, \quad i=1,\dots,n.
\ee
Below,  we reduce the problem of estimating vector function $\bof$ to estimating matrix $\bA_{0}$
of coefficients of $\bof$.

\subsection{Tensor approach to estimation}
\label{sec:tensor}

Denote $\ba = \Vect(\bA_{0})$ and   $\boB_i =  \Vect(\bphi(t_i) \bW_i^T)$. 
Note that $\boB_i$ is the $p(L+1)$-dimensional   vector with components
$\phi_l(t_i) \bW_i^{(j)}$, $l=0, \cdots, L$, $j=1, \cdots, p$, 
where $\bW_i^{(j)}$ is the $j$-th component of vector $\bW_i$. 
Consider matrix $\boB \in \RR^{n\times p(L+1)}$ with rows $\boB_i^T$, $i = 1, \cdots, n$,
vector $\bxi = (\xi_1, \cdots, \xi_n)^T$ 
and   vector $\bbb$ with components  $\bbb_i = \bW_i^T \brho (t_i)$, $i = 1, \cdots, n$.
Taking into account that 
$$
\Tr (\bA^T \bphi(t_i) \bW_i^T) =  \boB_i^T \Vect(\bA)
$$
we rewrite the varying coefficient model \fr{obs_mod} in a matrix form
\be \label{vcm_matr}
 \bY = \boB \ba + \bbb + \sig \bxi.
\ee
In what follows, we denote 
\be \label{matr_notations}
\bOmega_i = \bW_i \bW_i^T, \quad  
\bPhi_i = \bphi(t_i) (\bphi(t_i))^T, \quad 
\bSigma_i = \bOmega_i \otimes \bPhi_i,
\ee
where $\bOmega_i \otimes \bPhi_i$ is the Kronecker product of $\bOmega_i$ and $\bPhi_i$.
Note that $\bOmega_i$, $\bPhi_i$ and $\bSigma_i$ are i.i.d. for $i = 1, \cdots, n$, and that 
$\bOmega_{i_1}$ and $\bPhi_{i_2}$ are independent for any $i_1$ and $i_2$.  
By simple calculations, we derive  
\begin{eqnarray*}
\ba^T \boB \boB^T \ba &  = & \sum_{i=1}^n (\boB_i^T \ba)^2 
= \sum_{i=1}^n \lkv \Tr (\bA^T \bphi(t_i) \bW_i^T) \rkv^2 \\
& = & \sum_{i=1}^n \bW_i^T \bA^T \bphi(t_i)  \bphi^T (t_i)\bA \bW_i  =
\sum_{i=1}^n \ba^T (\bOmega_i \otimes \bPhi_i) \ba,
\end{eqnarray*}
which implies 
\be \label{eq:BtB}
\boB^{T}\boB=\sum_{i=1}^n \bOmega_i \otimes \bPhi_i.
\ee 
Let
\be \label{hatsigma}
\hbSig = n^{-1}\boB^{T}\boB  = 
n^{-1}   \sum_{i=1}^n  \bSigma_i.
\ee
Then, due to the i.i.d. structure of the observations, one has 
\be \label{sigma}
\bSigma = \EE \bSigma_1 = \bOmega  \otimes \bPhi \quad \mbox{with} \quad 
\bOmega =   \EE (\bW_1  \bW_1^T) \quad \mbox{and} \quad \bPhi = \EE (\bphi(t_1) \bphi^T(t_1)).
\ee

\section{Assumptions and notations}
\label{sec:assump}
\setcounter{equation}{0}

\subsection{Notations}

In what follows, we use  bold script  for matrices and vectors, e.g., $\bA$  or $\ba$, 
and superscripts to denote elements of those matrices and vectors, e.g., $\bA^{(i,j)}$  or $\ba^{(j)}$.
Below, we   provide a brief summary of the notation used throughout this paper. 
 
\begin{itemize}
\item 
For any vector  $\ba \in \RR^{p}$,   denote  the standard $l_1$ and $l_2$ vector norms by 
$\| \ba \|_{1}$ and $\| \ba \|_{2}$, respectively. For vectors $\ba, \bc \in \RR^{p}$, denote their scalar product by 
$\langle \ba, \bc \rangle$.
\item 
  For any function $q(t)$, $t \in [0,1]$,   $\| q \|_{2}$ and $\left\langle\cdot\,,\cdot\right\rangle_{2}$ 
are, respectively, the norm and the scalar product in the space $L_2 ([0,1])$. Also, $\| q \|_{\infty} = \sup_{t \in [0,1]} |q(t)|.$

\item 
For any vector function $\bq(t) = (q_1(t), \cdots, q_p(t))^T$, denote 
$$
\|\bq(t)\|_2 = \lkv \sum_{j=1}^p \|q_j\|_2^2 \rkv^{1/2}.
$$
\item 
For any matrix $\bA$, denote its spectral and Frobenius norms by 
$\| \bA\|$ and $\| \bA\|_2$, respectively. 
\item Denote  the $k\times k$ identity  matrix by $\II_{k}$.
\item  For any numbers, $a$ and $b$, denote $a \vee b = \max(a,b)$ and $a \wedge b = \min(a,b)$.
\item 
In what follows, we use the symbol $C$ for a generic positive
constant, which is independent of $n$, $p$, $s$ and $l$, and may take different
values at different places.
\item 
If $\br = (r_1, \cdots, r_p)^T$ and $r_j ' = r_j  +1/2-1/\nu_j$ for some $1 \leq \nu_j < \infty$, denote 
$r^{*}_j = r_j\wedge r'_j$ and $r_{\min}^{*} = \min_j r^{*}_j$.   

\item Denote
\be \label{t-and-W}
\bt = (t_1, \cdots, t_n), \quad \WW = (\bW_1, \cdots, \bW_n),
\ee
i.e., $\WW$ is the $p\times n$ matrix with columns $\bW_i$, $i=1, \cdots, n$. 

%
%
\end{itemize}

\subsection{Assumptions}

We impose the following assumptions on the varying coefficient model \fr{obs_mod}.
\\

\noindent
{\bf (A0). } Only  $s$ out of $p$ functions $f_j$ are non-constant and depend on the time variable 
$t$, 	$s_0$ functions are constant and independent of $t$ and $(p-s-s_0)$ functions are identically equal to zero.  
 We denote by $\mathcal J$ the set of indices corresponding to the non-constant functions   $f_j$.
\\

\noindent
{\bf (A1). } 
Functions $(\phi_k(\cdot))_{k=0,\dots,\infty}$ form  an orthonormal   basis of $L_2 ([0,1])$,
 and are such that  $\phi_0(t) = 1$ and, for any $t \in [0,1]$,  
any $l \geq 0$ and some $C_\phi < \infty$
\be \label{conA1}
\sum_{k=0}^l \phi_k^2 (t) \leq C_\phi^2 (l+1).
\ee
\\

\noindent
{\bf (A2). }  The probability density function $g(t)$  is bounded above and below 
$0< g_1 \leq g(t) \leq g_2 < \infty$.   Moreover, the eigenvalues of
$\EE (\bphi   \bphi^T) = \bPhi$  are    bounded from above and below
$$
0 < \phimin = \lam_{\min} (\bPhi) \leq \lam_{\max} (\bPhi) = \phimax  < \infty.
$$ 
Here, $\phimin$ and $\phimax$ are absolute constants independent of $L$.
\\

\noindent

\noindent
{\bf (A3). } Functions $f_j(t)$ have efficient representations in basis $\phi_l$, in particular, 
for any $j = 1, \cdots, p$,  one has
\be \label{conA4}
\sum_{k=0}^\infty |a_{jk}|^{\nu_j} (k+1)^{\nu_j r_j '} \leq (C_{a})^{\nu_j}, \quad r_j ' = r_j  +1/2-1/\nu_j,
\ee
 for some $C_a>0$, $1 \leq \nu_j < \infty$ and $r_j  > \min(1/2, 1/\nu_j)$.  
In particular, if function $f_j(t)$ is constant or vanishes, then $r_j = \infty$.
  We denote vectors with elements $\nu_j$ and $r_j$, $j=1, \cdots, p$,  by $\bnu$ and $\br$, respectively, and the set of indices 
of finite elements $r_j$ by $\calJ$:
\be \label{Ups_def}
\calJ = \{j:\ r_j < \infty\}.
\ee
\\

  \noindent
  {\bf (A4). } Variables $\xi_i$, $i=1, \cdots, n$, are i.i.d. \textit{sub-Gaussian} such that
  $$\EE \xi_i=0, \quad \EE \xi^2_i=1 \quad \text{and}\quad \Vert\xi_i\Vert_{\psi_2}\leq K$$
  where $\Vert\cdot\Vert_{\psi_2}$ denotes the sub-Gaussian norm. 
  \\

  \noindent
 {\bf (A5) }  ``Restricted Isometry in Expectation'' condition. Let $\bW_{\Lambda},\, \Lambda\in \{1,\dots,p\}$ 
be the sub-vector obtained by extracting the entries of $\bW$ corresponding  to indices in $\Lambda$ and let 
$\bOmega_{\Lambda}= \EE \left (\bW_{\Lambda}\bW_{\Lambda}^{T}\right )$.
   We assume that there exist two positive constants $\omax(\aleph)$ and $\omin(\aleph)$ such that for all subsets 
$\Lambda$ with cardinality  $\vert \Lambda\vert \leq \aleph$ and all  $\mathbf{v}\in \mathbf{R}^{|\Lambda|}$  one has
    \begin{equation}\label{RICE}
  \omin(\aleph)\Vert \mathbf{v}\Vert_{2}^{2}\leq  \mathbf{v}^T \bOmega_\Lambda \mathbf{v} \leq \omax(\aleph)\Vert \mathbf{v}\Vert_{2}^{2}. 
    \end{equation}
    Moreover, we suppose  that  $\EE (\bW^{(j)})^4 \leq V$ for any $j=1, \cdots, p$, 
    and that, for any $\mu\geq 1$ and for all subsets $\Lambda$ with $\vert \Lambda\vert \leq \aleph$, there exist positive constants 
        $U_\mu$ and $C_\mu$  and a set $\wpmu$ such that 
        \be \label{conA3pr}
        \bW_{\Lambda} \in \wpmu \Longrightarrow \lkr \| \bW_{\Lambda} \|_2 
        \leq U_{\mu} \rkr \cap \lkr  \max_{j \in \Lambda} |\bW_{\Lambda}^{(j)}|  \leq C_{\mu} \rkr,
        \quad \PP \lkr \bW_{\Lambda} \in \wpmu \rkr \geq 1-2\,p^{-2\mu}.
        \ee
         Here,     $\Umu = \Umu(\aleph)$, 
        $\Cmu= \Cmu(\aleph) $.\\ 
        \\
        
        \noindent
        {\bf (A6). }  We assume that $(s+s_0)(1+\log n)\leq p$ and that there exists a numerical constant $C_{\om}>1$ such that
        \begin{equation*}
        \log(n) \geq \dfrac{C_{\om}\,\phimax\,\omax((s+s_0)\log n)}{\phimin\ \omin((s+s_0)(1+\log n))}.
        \end{equation*}
        We denote
        \begin{equation}\label{ominmax}
         \omaxs =\omax\left ((s+s_0)(1+\log n)\right ), \quad\omins =\omin\left ((s+s_0)(1+\log n)\right ).
        \end{equation} 

\subsection{Discussion of assumptions}

\begin{itemize}
\item
Assumptions {\bf (A0)} corresponds to the case when $s$ of the covariates $f_j(t)$
are indeed functions of time, $s_0$ of them are time independent and $(p - s-s_0)$ are irrelevant.
\\ 

\item
Assumption {\bf (A1)}  deals with the basis of $L_2 ([0,1])$.
There are many types of orthonormal bases satisfying those conditions.
\\

\noindent
a) {\it Fourier basis}. Choose  $\phi_0(t) =1$, $\phi_k(t) = 2 \sin(2\pi k t)$ 
if $k>1$ is odd, $\phi_k(t) =  2 \cos(2\pi kt)$
if $k>1$ is even. The basis functions are bounded and     $C_\phi = 2$. 
\\

\noindent
b) {\it Wavelet basis}.  Consider a periodic wavelet basis on $[0,1]$:  
   $\psi_{h,i} (t)  = 2^{h/2} \psi(2^h t - i)$  
with $h=0, 1, \cdots$, $i=0, \cdots, 2^h-1$. Set  $\phi_0(t) =1$ and   
$\phi_j (t) =  \psi_{h, i} (t)$ with $j = 2^h + i +1$.   
If $l = 2^J$, then the condition \eqref{conA1} is satisfied with $C_\phi =  \| \psi \|_\infty$. Observing that, for $2^J < l < 2^{J+1}$, we have
$(l+1) \geq (2^{J+1} + 1)/2$  one can take  $C_\phi = 2 \| \psi \|_\infty$.
\\

\item Assumption {\bf (A2)} that $\phimin$ and $\phimax$ are absolute constants independent of $L$ is guaranteed by the fact that the  
sampling density $g$ is bounded above and below. For example, if $g(t) =1$, one has $\phimin = \phimax =1$.

\item
Assumption  {\bf (A3)}  describes sparsity of the vectors of coefficients of functions $f_j(t)$ in basis $\phi_l$, $j=1, \cdots, p$ and its smoothness. For example, 
if $\nu_j <2$, the vector of coefficients $a_{jl}$  of   $f_j$  is sparse.  
In the case when basis $\phi_l$ is formed by wavelets, condition    \fr{conA4}  implies that $f_j$ belongs to a Besov ball
of radius $C_a$. If we chose Fourier bases and $\nu_j = 2$, then   $f_j$ belongs to a Sobolev ball of smoothness  $r_j $ and radius $C_a$.
 Note that Assumption    {\bf (A3)} allows each non-constant function $f_j$ to have its own sparsity and smoothness patterns.

\item
Assumption {\bf (A4)} that $\xi_i$ are  sub-Gaussian random variables means that   
 their distribution  is dominated by the distribution of a centered Gaussian random variable. 
 This is a convenient and reasonably wide class. Classical examples of sub-gaussian random variables are Gaussian, Bernoulli
 and all bounded random variables. Note that $\EE \xi^2_i=1$ implies that $K\leq 1$.  

 \item 
  Assumption {\bf (A5)} is closely related to the Restricted Isometry (RI) conditions usually considered in the papers 
that employ LASSO technique or its versions, see, e.g., \cite{bickel_ritov_tsybakov}.  However, usually
the RI condition  is imposed on the matrix of scalar products of the elements of a deterministic  dictionary while we 
deal with a random dictionary and require this condition  to hold only for the expectation of this matrix.

Note that the upper bound in the condition \eqref{RICE} is automatically satisfied with $\omax=\Vert  \bOmega\Vert $ where
 $\Vert \bOmega\Vert$ is the spectral norm of the matrix $\EE (\bW   \bW^T) = \bOmega$. If the smallest eigenvalue of $\bOmega$,  $\lambda_{\min}(\bOmega)$, is non-zero, then the lower bound in \eqref{RICE} is satisfied with   $\omin=\lambda_{\min}(\bOmega)$.  
However, since  the $\lambda-$restricted maximal eigenvalue $\omax(\lambda)$ may be much smaller than the spectral norm of $\bOmega$ and 
$\omin(\lambda)$ may be much larger then $\lambda_{\min}(\bOmega)$, using those values will result in sharper  bounds for the error. 
Note that in the case when $\bW$ has i.i.d. zero-mean entries $W^{j}$ with $\EE \left (W^{(j)}\right )^{2}=\nu^{2}$, we have $\omax=\omin=\nu^{2}$.

 \item 
Assumption {\bf (A6)}   is usual in the literature on the hight-dimensional linear regression model, see, e.g., \cite{bickel_ritov_tsybakov}. 
For instance, if  $\bW$ has i.i.d. zero-mean entries $W^{j}$  and  $g(t) =1$, this condition is satisfied for any $1 < C_{\om} \leq \log(n)$.

\end{itemize}

\section{Estimation strategy and non-asymptotic error bounds }
\label{sec:est-err-bou}
\setcounter{equation}{0}

\subsection{Estimation strategy}
\label{sec:est-strategy}

Formulation \fr{vcm_matr} implies that the varying coefficient model
can be reduced to the linear regression model and one can apply one of the multitude  of 
penalized optimization techniques which have been developed for the linear regression. 
In what follows, we apply a block LASSO penalties for the   coefficients  in order to account 
for both the constant and the vanishing functions $f_j$
and also to take advantage of the sparsity of the functional coefficients in the chosen basis.

In particular, for each function $f_j$, $j=1, \cdots, p$, we divide its coefficients into $M+1$
different groups where group zero contains only coefficient $a_{j0}$ for the constant function $\phi_0(t) =1$ 
and $M$ groups of size $d \approx \log n$ where $M = L/d$. We denote    $\ba_{j0} = a_{j0}$ and
$\ba_{jl} = (a_{j,d(l-1)+1}, \cdots, a_{j,dl})^T$ the sub-vector of coefficients  of function $f_j$ in block $l$, $l=1, \cdots, M$.
Let $K_l$ be the subset of indices associated with  $\ba_{jl}$.
 We impose block norm on matrix $\bA$  as follows
\be \label{block}
\| \bA \|_{\block} = \sum_{j=1}^p \sum_{l=0}^M \| \ba_{jl} \|_2.
\ee
Observe that $\| \bA \|_{\block}$   indeed satisfies the definition of a norm and is a sum  of absolute values of  
 coefficients $a_{j0}$  of functions $f_j$  and $l_2$ norms for each of the block vectors 
of coefficients $\ba_{jl}$, $j=1, \cdots, p$, $l=1, \cdots, M$. 


The penalty which we impose is related 
to both the ordinary and the group LASSO penalties which have been used by many authors. 
The difference, however, lies in the fact that the structure of 
the blocks is not motivated by naturally occurring groups   (like, e.g., rows of the matrix $\bA$) but rather our desire to 
exploit sparsity of functional coefficients $a_{jl}$. In particular, we construct an estimator $\hbA$ of $\bA_{0}$ as a solution 
of the following convex optimization problem
\be \label{matr_opt}
\hbA = \arg\min_{\bA} \lfi n^{-1} \sum_{i=1}^n  \lkv Y_i - \Tr(\bA^T \bphi(t_i) \bW_i^T) \rkv^2 + \delta \| \bA \|_{\block} \rfi,
\ee
where the value of $\delta$ will be defined later.

Note that with the tensor approach which we used in Section \ref{sec:tensor}, optimization problem 
\fr{matr_opt} can be re-written in terms of vector    $\balpha = \Vect(\bA)$   as
 \be \label{vect_opt}
\hba = \arg\min_{\balpha} \lfi n^{-1} \| \bY - \boB \balpha \|_2^2 + \del \| \balpha \|_{block} \rfi,
\ee
where $\| \balpha \|_{block} = \| \bA \|_{block}$ is defined by the right-hand side  \fr{block} with vectors 
$\ba_{jl}$ being sub-vectors of vector $\balpha$.  Subsequently, we construct an estimator  
$\hbf (t) = (\hf_1 (t), \cdots, \hf_p(t))^T$ of the vector function $\bof(t)$ using
\be \label{fun_coef}
\hf_j(t) = \sum_{k=0}^L \ha_{jk} \phi_k (t), \quad j=1, \cdots, p.
\ee

In what follows, we derive the upper bounds for the risk of the estimator  $\hba$ (or $\hbA$) and suggest a value 
of parameter $\del$ which allows to attain those bounds. However, in order to obtain a benchmark of how 
well the procedure is performing, we determine the  lower bounds for the risk of any estimator $\hbA$ under assumptions
{\bf (A0)}--{\bf (A4)}.

\begin{remark} \label{remarkL}
{\rm Assumption that   $K = L/d$  is an integer  is not essential. Indeed, we can replace  the number of groups $K$ 
by the largest integer below or equal to $L/d$ and then adjust group sizes to be $d$ or $d+1$
where $d =[\log n]$, the largest integer not exceeding $\log n$. }
\end{remark}

\subsection{Lower bounds for the risk}
\label{se:low_bou}
In this section we will obtain the lower bounds on the estimation risk. 
We consider a class ${\cal F} = {\cal F}_{s_0,s, \bnu, \br } (C_a)$ of vector functions $\bof  (t)$ such that $s$ 
of their components  are non-constant with coefficients satisfying condition \fr{conA4} in {\bf (A3)}, 	
$s_0$ of the components are constant and  $(p-s-s_0)$ components are identically equal to zero.
We construct the lower bound  for the minimax quadratic risk of any estimator  $\tbf$
of the vector function $\bof \in {\cal F}_{s_0, s,\bnu, \br } (C_a)$.
 Let $\omaxs$ be given by formula \fr{ominmax}. Denote $\rmax = \max \{r_j:  j \in \calJ \}$  and
\be \label{n_low}
n_{low} = \frac{  2 \, \sig^2 \kappa}{C_a^2\, \omaxs \,\phimax} \,  \max\left \{1,\lkr \frac{6}{s} \rkr^{2 \rmax  +1}\right \},
\ee
 \be \label{Delta_low}
\Delta_{lower} (s_0, s, n, \br ) = \max \lfi \frac{\kappa\,\sigma^2\,  s_0   }{ 4  n\, \omaxs\, \phimax},\  
  \frac{1}{8}\ 
\sum_{j \in \calJ} \ C_a^{\frac{2}{2 r_j   +1}}\ \lkr \frac{\sigma^2\,\kappa   }{n \,\omaxs \,\phimax} \rkr^{\frac{2r_j }{2 r_j +1}}    \rfi.
\ee 
Then, the following statement holds.

\begin{theorem}  \label{th:lowbound}
 Let   $s\geq 1$ and $s_0\geq 3$.  
Consider observations $Y_i$  in model \fr{obs_mod} with   $\bW_i$,  $i=1, \cdots, n$ and $t$  satisfying assumptions  
{\bf (A5)}  and {\bf (A2)}, respectively.   Assume that, conditionally on $\bW_i$ and $t_i$,   variables $\xi_i$ are Gaussian $\mathcal{N}(0,1)$ and 
that $n\geq n_{low}$. Then, for  any  $\kappa < 1/8$ and   any  estimator  $\tbf$  of   $\bof$, one has   
\be \label{eq:loweb}
  \inf_{\tbf} \, \sup_{\bof \in  {\cal F}  } 
\PP \lkr \|\tbf - \bof\|_2^2 \geq \Delta_{lower} (s_0,  s, n, r ) \rkr \geq 
 \frac{\sqrt 2}{1 + \sqrt 2} 
 \lkr 1 - 2\kappa - \sqrt{\frac{2 \kappa}{\log 2}} \rkr.   
\ee
\end{theorem}

 Note that condition  $ s_0 \geq 3$  is not essential since, for $s_0  < 3$, the first term in \fr{Delta_low} is of parametric 
order. Condition  $n \geq n_{low}$ is a  purely technical condition which is satisfied for the collection of $n$'s 
for which upper bounds are derived.
Observe also that inequality \fr{eq:loweb} immediately implies that 
\be \label{eq:lowmean}
 \inf_{\tbf} \, \sup_{\bof \in  {\cal F}  } 
\EE  \|\tbf - \bof\|_2^2 \geq \Delta_{lower} (s_0,  s, n, r )  
 \lkv  
 \frac{\sqrt 2}{1 + \sqrt 2} 
 \lkr 1 - 2\kappa - \sqrt{\frac{2 \kappa}{\log 2}} \rkr \rkv.   
\ee

\subsection{ Adaptive estimation and upper bounds for the risk } 
\label{se:up_bou}

In this section we derive an  upper bound for the risk of the estimator  \eqref{matr_opt}. 
For this purpose, first, we shall show that, with high probability,
   the ratio between the restricted eigenvalues of matrices $\hbSig$ defined in \fr{hatsigma} and $\bSigma = \EE \hbSig$
   is bounded above and below. 
%
    This is accomplished by the following lemma. 
    
     For any $\Lambda\in \{1,\dots,p\}$, we denote by
     $\left (\bOmega_{i}\right )_{\Lambda}=(\bW_{i})_{\Lambda}\left ((\bW_{i})_{\Lambda}\right )^{T}$,  
    $\hbSig_{\Lambda} =n^{-1}\sum_{i=1}^n (\bOmega_i )_{\Lambda}\otimes \bPhi_i$ and $\bSig_{\Lambda}=\EE \hbSig_{\Lambda} $ .   
     For  some   $0 < h < 1$ and $1\leq\aleph\leq p$, we define 
       \be \label{No}
       N(\aleph) = \frac{64\,\mu\,\aleph\,  (L+1)\, \log(p+L)\ C_\phi^2\, U_\mu^2(\aleph)\, \phimax\, \omax(\aleph)  }
       {h^2\,  \phimin^2\, \omin^{2}(\aleph)}.
       \ee

    \begin{lemma} \label{lem:Sigma} 
   Let $n  \geq N(\aleph)$ and $\mu$ in \fr{conA3pr} be large enough, so that 
   \be \label{mu_cond}
   p^{\mu} \geq \max\left \{\frac{\sqrt{2V}\,n} {8\mu\sqrt{\aleph}\,U^{2}_{\mu} (\aleph)\log(p+L)}, \  2n\right \},
   \ee
  where $V$ is defined in Assumption (A5).
   Then,  for any $\Lambda\in \{1,\dots,p\}$ 
   \be \label{LargeDevSig} 
   \inf_{ \Lam:\, | \Lam |  \leq \aleph }\   
   \left[ \PP \lkr \lfi  \| \hbSig_{\Lambda} - \bSig_{\Lambda} \| < h\, \phimin\, \omin(\aleph) \rfi \cap  \wpmun  \rkr \right]
   \geq 1 -2p^{-\mu}, 
   \ee
    where   $\wpmu$ is the set of points in $\RR^p$ such that condition \fr{conA3pr} holds
      and   $\wpmun$ is  the direct product of $n$ sets $\wpmu$. 
      
      Moreover, on the  set $\wpmun$, with probability at least $1 -2p^{-\mu}$, one has simultaneously
   \be \label{eigenSig}
   \inf_{ \Lam:\, | \Lam |  \leq \aleph }\  \left[ \lamin(\hbSig_{\Lambda}) \right] \geq (1-h) \phimin\, \omin(\aleph), \quad
   \sup_{ \Lam:\, | \Lam |  \leq \aleph }\  \left[ \lamax(\hbSig_{\Lambda}) \right]  \leq (1+h) \phimax\, \omax(\aleph).
   \ee
    \end{lemma}

\medskip

Lemma  \ref{lem:Sigma} ensures that the restricted lowest eigenvalue of the regression matrix $\hbSig$ 
is within a constant factor of the respective eigenvalue of matrix $\bSig$. Since $p$ may be large,
this is not guaranteed  by a large value of $n$ (as it happens in the asymptotic setup) and leads to additional conditions  on the 
relationship between parameters $L$, $p$ and $n$. 
 
By applying a combination 
of LASSO and group LASSO arguments, we obtain the following theorem that gives an upper bound for the quadratic risk of the estimator \eqref{matr_opt}.
We set    $U_{\mu}=U_{\mu}(s+s_0)$. Define
\be \label{Nmax}  
\mathbf N = 
\max\left \{\frac{64\,\mu\,(s+s_0) C_\phi^2 U_\mu^2 (L+1) \phimax \omaxs \log(p+L) }{h^2 \phimin^2 (\omins)^{2}}, 
\dfrac{U_{\mu}^{2}\,C_{\phi}^{2}(L+1)\mu\log p}{g_2\, \omax (s)}, 3\,C^{2}_a\,g_2\,s\, \omax(s) \right \}
\ee
and 
\be \label{deldef}
\hdel = 2\,\left (\sig C_{\om} \,K\sqrt{\mu}+1)\right )  \sqrt{\frac{(1+h)\phimax\omax(1)\log p}{n}}.
\ee

\begin{theorem} \label{th:upper_bou}
 Let $\min_k (r_k \wedge r_k' )\geq 2$, $L+1 \geq  n^{1/2}$ and $n  \geq \mathbf N$.  Let $\mu$ in \fr{conA3pr} be large enough, so that 
   \be \label{mu_cond_2}
   p^{\mu} \geq \max\left \{\frac{\sqrt{2V}\,n} {8\mu\sqrt{s+s_0}\,U^{2}_{\mu} \log(p+L)},  \dfrac{2L}{\log n},2n \right \}.
   \ee
 If  $\hba$ is an estimator of $\ba$ obtained as a solution of optimization problem \fr{vect_opt} 
with $\del = \hdel$, and the vector function
$\hat{\mathbf f}$  is recovered using \fr{fun_coef}, then, one has 
\be \label{up_bou1} 
\PP \lkr  \| \hbf - \bof \|^2_2   \leq \Delta  (s_0, s, n, \br ) \rkr \geq   1 -   8\, p^{-\mu}
\ee 
where  
\begin{eqnarray*}  
\Delta  (s_0, s, n, \br ) & = & \dfrac{C^{2}_a\,s}{n^{2}}+
\frac{C_{B}\, (1+h) \omaxs \phimax}{(1-h) \omins \phimin} 
\lkv \frac{  \left (C_{\om}\, \sig^2  K^2\mu+1\right )  (s_0+s)\log p }{n (1-h) \omins \phimin} 
 \right. \\
& + & \left. \sum_{j \in \calJ}\   C_a^{2/(2r_j+1)}\   \lkr \frac {C_{\om}\, \sig^2  K^2\mu+1 }
{n (1-h) \omins \phimin}  \rkr^{\frac{2r_j}{2r_j +1}}
(\log n)^\frac{(2-\nu_j)_+-2\nu_j\,r_j}{\nu_j(2r_j +1)}\left (\log p\right )^{\frac{2r_j}{2r_j +1}} \rkv. 
\end{eqnarray*}
\end{theorem}

Note that  construction \fr{vect_opt}  of the estimator $\hba$ does not involve knowledge of   unknown parameters  
$\br$ and  $\bnu$  or matrix $\bSig$, therefore, estimator $\hba$ is fully adaptive. 
Moreover, conclusions of Theorem \ref{th:upper_bou} are derived without any asymptotic assumptions on $n$, $p$ and $L$.
\\

In order to assess the optimality of estimator $\hba$, we  consider the case of the Gaussian noise, i.e. $K=1$. 
Observe that, under Assumption {\bf (A2)}, 
the values of  $\phimin$ and $\phimax$ are independent of $n$ and $p$, so that the only quantities in \fr{up_bou1}  
which are not bounded from above and below by an absolute constant are $\sigma$, $\omins$, $\omaxs$, $s$ and $s_0$.  Hence, 
$\Delta  (s_0, s, n, \br ) \leq C\, \Delta_{upper}  (s_0, s, n, \br )$ with
\be \label{Delta_up1}
\Delta_{upper}  (s_0, s, n, \br )  =  \frac{\omaxs}{\omins} \lkv 
 \frac{\sig^2  s_0 \log p}{n \omins} +   
\sum_{j \in \calJ}\  C_a^{2/(2r_j+1)} \lkr \frac {\sig^2} {n  \omins}  \rkr^{\frac{2r_j}{2r_j +1}}\   
\left (\log p\right )^{\frac{2r_j}{2r_j +1}} 
 \rkv,
\ee
where $C$ is an absolute constant independent of $n$, $p$ and  $\sig^2$ and we use $(2-\nu_j)_+-2\nu_j\,r_j\leq 0$.

 Inequality \fr{Delta_up1} implies that, for any values of the parameters, the ratio  between the upper and the lower 
bound for the risk \fr{Delta_low} is bounded by $C \log( p) \, {\omaxs}^2/{\omins}^2$.
Note that $\omaxs/\omins$ is the condition number of matrix $\bOmega_{\Lam}$ with $|\Lam| = (s+s_0)(1+\log n)$. 
Hence, if matrix $\bOmega_{\Lam}$ is well conditioned,
so that $\omaxs/\omins$ is bounded by a constant,   the estimator $\hbf$ attains  optimal convergence rates up to a $\log p$ factor.

Consider the case when all functions $f_j(t)$ are spatially homogeneous, i.e., $\min_j \nu_j \geq 2$ and  $n$ is large enough, 
i. e.  there exists a positive $\beta$ such that $n^{\beta}\geq p$. Then, the estimator $\hbf$ attains  optimal convergence rates up to a
constant factor, if $s/s_0$ is bounded or $n \omins/\sig^2$ is relatively large. 
In particular,  if all functions in assumption {\bf (A3)} belong to the same space, then the following corollary is valid.

\begin{corollary} \label{cor:rates}
Let conditions of Theorem \ref{th:upper_bou} hold with $r_j =r$ and $\nu_j = \nu$, $j-1, \cdots, p$,
and   matrix $\bOmega$ be well conditioned, i.e. 
$\omaxs/\omins$ is bounded by for some absolute constant independent of $n$, $p$ and $\sigma$. Then, 
\be \label{error_ratio}
\frac{\Delta_{upper}  (s_0, s, n, \br)} {\Delta_{lower}  (s_0, s, n, \br )} \leq \lfi
\begin{array}{ll}
C   \log p, & \mbox{if} \quad  \sig^{2} (s/s_0)^{2r+1} \geq n \omins, \\
C    (\log p)^\frac{2r}{2r+1}, & \mbox{if} \quad  \sig^{2} (s/s_0)^{2r+1} <   n \omins, \ 1 \leq \nu <2,\\
C, & \mbox{if} \quad  \sig^{2} (s/s_0)^{2r+1} <   n \omins, \ \nu \geq 2.
\end{array} \right.
\ee
\end{corollary}

\subsection{Adaptive estimation with respect to the mean squared risk} 
\label{sec:mean_squared_risk}

Theorem \ref{th:upper_bou} derives upper bounds for the risk with high probability.
Suppose that an upper bound on the norms of functions  $\bof_j$ is available due to physical or other considerations:
\be \label{norm_upper}
\max_{1 \leq j \leq p} \|f_j \|_2^2 \leq C_f^2.
\ee
Then, $\| \ba \|^2 \leq p\,  C_f^2 $ and $\hba$ given by \fr{vect_opt} can be replaced by the solution of the convex problem
 \be \label{vect_opt_bounded}
 \hba = \arg\min_{\ba}  \lfi    n^{-1} \| \bY - \boB \ba \|_2^2 + \del \| \ba \|_{block}  \quad  \mbox{s.t.} \quad  \|\ba \|^2 \leq p \,C_f^2  \rfi  
\ee
 with $\del = \hdel$ where  $\hdel$  is defined in \fr{deldef}, 
and estimators $\hf_j$ of $f_j$, $j=1, \cdots, p$, are constructed using formula \fr{fun_coef}.
Choose $\mu$ in  \fr{conA3pr} large enough, so that 
\be \label{tau_bound}
16 n C_f^2 \leq p^{\mu -1}.
\ee
Then, the following statement is valid.

\begin{theorem} \label{th:upper_mean}
Under the assumptions of the Theorem \ref{th:upper_bou},
 and  for 
any $\mu$ satisfying condition \fr{tau_bound}, one has 
\be \label{up_bou_mean} 
\EE  \| \hbf - \bof \|^2_2   \leq C \Delta_{upper}  (s_0, s, n, \br )
\ee 
where $C$ is an absolute constant  independent of $n$, $p$ and $\sigma$.
\end{theorem}
 
\medskip


\section{Examples and discussion }
\label{sec:exam_discussion}
\setcounter{equation}{0}

\subsection{Examples}
\label{sec:examples}

In this section we provide several examples when assumptions of the paper are satisfied. 
For simplicity, we assume that $g(t) =1$, so that $\phimin = \phimax =1$.

\begin{example} \label{ex:normal} {\bf Normally distributed dictionary}
{\rm   Let $\bW_i$, $i=1, \cdots, n$, be i.i.d. standard Gaussian vectors $N(\bzero, \bI_p)$.
Then, $\bOmega = \bI_p$, so that $\omin = \omax = 1$. Moreover, $ \bW^{(j)}$ are independent 
standard Gaussian variables and $ (\bW_i)_{\Lambda}^T (\bW_i)_{\Lambda}$ are independent
chi-squared variables with $\vert \Lambda \vert=\aleph$ degrees of freedom.   
Using inequality  (see, e.g., \cite{birge}, page 67)
$$
\PP \lkr \chi_\aleph^2 \leq \aleph + 2 \sqrt{ \aleph x} + 2 x \rkr \geq 1 - e^{-x},\quad x>0,
$$
 for any  $\mu_1 \geq 0$, derive   
$$
\PP \lkr (\bW_1)_{\Lambda}^T (\bW_1)_{\Lambda} \leq (\sqrt{\aleph} + \sqrt{2} \mu_1)^2 \rkr \geq 1 - \exp(-  \mu_1^2).
$$
Choose any $\mu>0$ and set $\mu_1^2 =  2 \mu \log( p)$.
Then, using a standard bound on the maximum of $p$ Gaussian variables
one obtains that Assumption {\bf (A5)} holds with 
$$
C_{\mu} = \sqrt{ 2 \log p}, \quad  U_{\mu}^2  = (\sqrt{\aleph} + 2\sqrt{\mu \log p})^2.
$$} 
\end{example}

\begin{example} \label{ex:Bernoulli}{\bf Symmetric Bernoulli dictionary}
{\rm   Let $\bW_i^{(j)}$, $i=1, \cdots, n$, $j=1, \cdots, p$, be 
independent symmetric Bernoulli variables 
$$
\PP(\bW_i^{(j)} = 1) = \PP(\bW_i^{(j)} = - 1)=1/2.
$$
Then,  $\bOmega = \bI_p$,   $\omin = \omax = 1$ and, for any  $\mu$,  
$$C_{\mu} = 1, \quad  U_{\mu}^2(\aleph) = \aleph.
$$
 
} 
\end{example}

In both cases, $\mathbf N$ in \eqref{Nmax} is of the form $\mathbf N=C(s+s_0)^{2}(L+1)\log(p)$. 
Under conditions of Theorem \ref{th:upper_bou}, the   upper bounds for the risk are of the form
$$
 \Delta_{upper}  (s_0, s, n, \br ) = C 
\lkv  \frac{\sig^2  s_0 \log n}{n } +   
\sum_{j \in \calJ}\  \lkr \frac {\sig^2} {n}  \rkr^{\frac{2r_j}{2r_j +1}}\,C_a^{2/(2r_j+1)}   
(\log n)^\frac{(2-\nu_j)_+-2\nu_j\,r_j}{\nu_j(2r_j +1)}\left (\log p\right )^{\frac{2r_j}{2r_j +1}}  \rkv
$$
 where $C$ is a numerical constant, so that it follows from Corollary \ref{cor:rates}
that the  block  LASSO estimator is   minimax   optimal up to, at most logarithmic factor of $p$.

 The two examples above illustrate the situation when estimator \fr{fun_coef} attains nearly optimal convergence rates 
when $p>n$. This, however, is not always possible. Note that  our analysis of the performance of 
the estimator \fr{fun_coef} relies on the fact that the  eigenvalues of any sub-matrix $\hbSig_{\Lambda}$ 
are close to those of matrix $\bSig_{\Lambda}$  (Lemma \ref{lem:Sigma}). The latter requires $n \geq   \mathbf N$  
where  $\mathbf N$  depends on the nature of vectors  $\bW_i$. The next example shows that sometimes $n<p$ 
prevents Lemma \ref{lem:Sigma}  from being valid.

\begin{example} \label{ex:orthonorm} {\bf Orthonormal dictionary}
{\rm   Let $\bW_i$, $i=1, \cdots, n$, be uniformly distributed on a set of 
canonical vectors $\bone_k$, $k=1, \cdots, p$. Then, $\bOmega = \bI_p/p$, so that $\omin = \omax = 1/p$.
Moreover, $\| \bW_1\|_2^2 =1$ and $|\bW^{(j)}| \leq 1$.
Therefore, for any  $\mu>0$,  
$$
C_{\mu} = 1, \quad  U_{\mu}^2 (\aleph)= 1.
$$
\\

In the case of the orthonormal dictionary,   $\mathbf N$ in \eqref{Nmax} is of the form $\mathbf N=C(s+s_0)(L+1)p\log(p)$.
Under conditions of Theorem \ref{th:upper_bou}, the   upper bound for the risk is of the form
$$
 \Delta_{upper}  (s_0, s, n, \br ) = C 
\lkv  \frac{\sig^2 p s_0 \log n}{n} +   
\sum_{j \in \calJ}\  \lkr \frac {p \sig^2} {n}  \rkr^{\frac{2r_j}{2r_j +1}}\,C_a^{2/(2r_j+1)}   
(\log n)^\frac{(2-\nu_j)_+-2\nu_j\,r_j}{\nu_j(2r_j +1)}\left (\log p\right )^{\frac{2r_j}{2r_j +1}}\rkv,
$$ 
so,   $n\geq \mathbf N $ implies $n > C \sig^2 p (s_0+s) \log n$ which also guarantees that the risk of the estimator is small. 
This, indeed, coincides with one's intuition since one would need to sample more than $p$ vectors in order to ensure 
that each component of the vector has been sampled at least once.}
\end{example}
  
\subsection{Discussion}
\label{sec:discussion}

In the present paper, we provided a non-asymptotic minimax study of the sparse high-dimensional varying coefficient model. 
To the best of our knowledge, this has never been accomplished before.   An important feature  of our analysis is its flexibility:
it distinguishes between vanishing, constant and time-varying covariates and, in addition, it allows the latter to be heterogeneous 
(i.e., to have different degrees of smoothness) and spatially inhomogeneous.
In this sense, our setup is more flexible than the one usually  used in the 
context of additive or compound functional models  (see, e.g., \cite{dalalyan} or \cite{raskutti}).

The adaptive estimator is obtained  using block LASSO approach which can be viewed as a version of group LASSO 
where groups do not occur naturally but are rather driven by the need to reduce the variance, 
as it is done, for example,  in block thresholding.  Since we used tensor approach for derivation of the estimator, 
we believe that the results of the paper can be generalized to the case of the multivariate varying coefficient model 
studied in \cite{zhu}. 

  An important feature of our estimator is that it is fully adaptive. Indeed, application of the  proposed
block LASSO technique does not require the knowledge of the number  of the non-zero components of $\bof$.  It only depends 
on the highest diagonal element of matrix $\bOmega$  which can be estimated with high precision even when $n$
is quite small due to Lemma \ref{lem:Sigma}.

  Note that, even when   $p$ is larger than $n$, the vector function  $\bof$ is completely 
identifiable due to Assumption {\bf (A5)}, as long as the number of non-zero components 
of  $\bof$ does not exceed $\aleph$ in condition \fr{RICE} and the number of observations $n$ is large enough. 
The examples of the paper deal with the dictionaries such that \fr{RICE} holds for all $\aleph =1, \cdots, p$.
The latter ensures identifiability of $\bof$ provided $n \geq  \mathbf N$ where $\mathbf N$  
is specified for each type of the random dictionary and depends on the sparsity level of $\bof$.
On the other hand, large values of $p$ ensure  great flexibility of the choice of $\bof$, 
so one can hope to represent the data   using only few components of it.

Finally, we want to comment on the situation when the requirement $n \geq  \mathbf N$ is not met
due to lack of sparsity or insufficient number of observations.
In this case,   $\bof$  is not identifiable and one cannot guarantee that $\hbf$ is close to 
the true function $\bof$. However, this kind of situations   occur  in all types of high-dimensional problems.

\section*{Acknowledgments}

Marianna Pensky  was partially supported by National Science Foundation
(NSF), grant  DMS-1106564. The authors want to thank Alexandre Tsybakov for 
extremely valuable suggestions and discussions.


\section{Proofs}
\label{sec:proofs}
\setcounter{equation}{0}


\subsection{  Proofs of the lower bounds for the risk}
\label{sec:proof-low-bou}

In order to prove Theorem \ref{th:lowbound}, we consider a set of test vector functions
$\bof_{\bom} (t) = (f_{1, \bom}, \cdots, f_{p, \bom})^T $
 indexed by binary sequences $\bom$ with components
\be \label{test-func}
f_{k, \bom} (t) = \om_{k0} u_k  + \sum_{l = l_{0k}}^{2 l_{0k}-1}  \om_{kl} v_k  \phi_l(t),   
\ee
where  $\om_{kl}  \in  \{0,1\}$ for  $l= l_{0k}, l_{0k} +1, \cdots, 2 l_{0k}-1$, $k=1, \cdots, p$. 
Let  $K_0$ and $K_1$, respectively, be the sets of indices such that 
$K_0 \cap K_1 = \emptyset$ and  $u_k = u$ if $k \in K_1$ and $u_k =0$ otherwise,
$v_k = v$ if $k \in K_0$ and $v_k =0$ otherwise, so that $u_k v_k=0$.

In order assumption \fr{conA4} holds, one needs $u  \leq C_a$ and 
\be \label{v_condition}
 v^{\nuj}   \sum_{l = l_{0k}}^{2 l_{0k}-1} (l+1)^{\nuj \rj '} \leq (C_a)^\nuj, \ j \in \Upsilon.
\ee
By simple calculations, it is easy to verify that   condition \fr{v_condition} is satisfied if we set 
\be \label{vk-coef}
u   \leq   C_a, \quad  v = C_a (2 l_{0k})^{-(r_k +1/2)},  
\ee
where the constancy of $v$ implies that   $l_{0k}$ in \fr{vk-coef} are different for different values of $k$.

 Consider two binary sequences $\bom$ and $\tbom$ and the corresponding test functions $\bof   (t) = \bof_{\bom} (t)$
and  $\tbf (t) = \bof_{\tbom} (t)$  indexed by those sequences. 
Then, the total squared distance   in $L_2 ([0,1])$  between $\bof_{\bom} (t)$ and  $\bof_{\tbom} (t)$ is equal to
\be \label{total_dist}
D^2 = u^2\ \sum_{k \in K_1}  |\om_{k0} - \tom_{k0}| +
v^2\ \sum_{k \in K_0} \sum_{l = l_{0k}}^{2 l_{0k}-1}  |\om_{kl} - \tom_{kl}|.  
\ee

  Let  $P_{\bof}$ and $P_{\tbof}$ be probability measures corresponding to test functions $\bof$ and $\tbof$, respectively. 
  We shall consider two cases. In    {\it Case 1}, the first $s_0$ functions are constant and  the rest of the functions are equal to identical zero. 
In   {\it Case 2}, the first $s$ functions are time-dependent and the rest of the functions are equal to identical zero. 
In both cases, $\bof   (t)$ and $\tbof (t)$ contain at most $s+s_0$ non-zero coordinates.
Using  that conditionally on $W_i$ and $t_i$, the variables $\xi_i$ are Gaussian $\mathcal{N}(0,1)$, we obtain  that the 
Kullback-Leibler divergence ${\cal K}(P_{\bof}, P_{\tbof})$ between $P_{\bof}$ and $P_{\tbof}$ satisfies
$$
{\cal K}(P_{\bof}, P_{\tbof}) = (2 \sig^2)^{-1}\, \EE  \sum_{i=1}^n \lkv Q_i (\bof) - Q_i (\tbof) \rkv^2 =
(2 \sig^2)^{-1} n\ \EE \lkv Q_1 (\bof) - Q_1 (\tbof) \rkv^2
$$
where 
$$
Q_i (\bof) =  W_i^{T}\bof(t_i)  
$$ 
and,  due to conditions {\bf (A2)} and {\bf (A5)},
 \begin{eqnarray*}
 \EE \lkv Q_1 (\bof) - Q_1 (\tbof) \rkv^2 & = & \EE \left ((\bof-\tbof)^{T}(t_1)W_1 W_1^{T}(\bof-\tbof)(t_1)\right ) 
=\EE \left ((\bof-\tbof)^{T}(t_1)\,\bOmega\, (\bof-\tbof)(t_1)\right ) \\
 & \leq & \omaxs\  \EE \left ( \left\Vert (\bof-\tbof)(t_1) \right\Vert^{2}_{2}\right )\leq  \omaxs\, \phimax D^2,
 \end{eqnarray*}
where $\omaxs$ and $D^2$ are defined in \fr{ominmax} and \fr{total_dist}, respectively.

In order to derive the lower bounds for the risk, we  use Theorem 2.5 of Tsybakov (2009) 
which implies that, if a set $\cOm$ of cardinality $M+1$ contains   sequences $\bom_0, \cdots, \bom_M$ with $M \geq 2$
such that, for any $j = 1, \cdots, M$, one has $\|f_{\bom_0} - f_{\bom_j}| \geq D >0$,  $P_{\bom_j} << P_{\bom_0}$ 
and 
${\cal K}(P_{\bof_j}, P_{\bof_0}) \leq \kappa \log M$ with $0 < \kappa < 1/8$,
then 
\be \label{tsybakov}
 \inf_{\tbom} \, \sup_{f_{\bom}, \bom \in \cOm}\ 
\PP \lkr \|\bof_{\bom} -\bof_{\tbom}\|_2  \geq D/2 \rkr \geq \frac{\sqrt M}{1 + \sqrt M} 
\lkr 1 - 2\kappa - \sqrt{\frac{2 \kappa}{\log M}} \rkr.   
\ee
Now, we consider two separate cases.
\\
 
{\it Case 1}. Let the first $s_0$ functions be constant and  the rest of the functions be equal to identical zero.
Then, $v=0$ and $K_1 = \{1, \cdots,  s_0 \}$.
Use the Varshamov-Gilbert Lemma (Lemma 2.9 of \cite{tsybakov_book}) to choose a set  $\cOm$ of $\omega$
 with  $\card(\cOm) \geq 2^{ s_0 /8}$ and $D^2 \geq u^2  s_0/8.$
Inequality
 
$$
{\cal K}(P_{\bof_j}, P_{\bof_0}) \leq (2 \sig^2)^{-1} n\,  
\omaxs\, \phimax u^2\,  s_0 /8 \leq \kappa\ \log  \left ( \card(\cOm)\right )
$$ 
holds if  $u^2  = 2\,\sig^2 \kappa /(n\, \omaxs\, \phimax)$. 
Then, $D^2 = (s_0 \sig^2 \kappa) /(4 n\, \omaxs \,\phimax)$ and $u  \leq C_a$  provided $n\geq 2\,\sig^2 \kappa /(C_a^{2}\omaxs \phimax)$.  
\\

{\it Case 2}. Let the first $s$ functions be time-dependent and the rest of the functions be equal to identical zero.
Then $u=0$, $v$ is given by formula \fr{vk-coef}, $K_0 = \{1, \cdots, s\}$. Let $r_k, k \in K_0$, coincide with the 
values of finite components of vector $\br$. Denote 
$$
\calL = \sum_{k=1}^s l_{0k}.
$$
Use Varshamov-Gilbert Lemma to choose a set $\cOm$ of $\omega$
 with  $\card(\cOm) \geq 2^{\calL /8}$ and $D^2 \geq v^2 \calL/8.$
Inequality ${\cal K}(P_{\bof_j}, P_{\bof_0}) \leq     \kappa \calL /8$  holds if 
$$
v^2   \leq ( \sig^2 \kappa) /(4 n\, \omaxs \,\phimax),
$$
which, together with \fr{vk-coef} and \fr{total_dist} imply that
$$
l_{0k}  =  \left  \lfloor\frac{1}{2} \lkr \frac{ 4\,C_a^2\, n\,  \omaxs\, \phimax}{\sig^2 \kappa} \rkr^{\frac{1}{2 r_k  +1}} \right \rfloor +1,
\quad
D^2 \geq \frac{C_a^2}{16}\ \sum_{k=1}^s \lkr \frac{4\,C_a^2\, n\, \omaxs\, \phimax}{\sig^2 \kappa} \rkr^{-\frac{2 r_k }{2 r_k  +1}} 
$$
where $\lfloor x\rfloor$ denotes the integer part of $x$. 
Condition  $\calL \geq 3$ is satisfied for any $s\geq 1$  provided  $n \geq n_{low}$.



\subsection{Proofs of the upper bounds for the risk}
\label{sec:proofs-upper}

{\bf Proof of Theorem \ref{th:upper_bou}. }\\

\noindent
 For any $\balpha \in \RR^{p(L+1)}$, one has  
$$
n^{-1} \| \boB \hba - \bY \|_2^2 +   \del \| \hba  \|_{block}  \leq 
n^{-1} \| \boB  \balpha  - \bY \|_2^2 +    \del \| \balpha   \|_{block}. 
$$
 Consider a set  $\calF_1 \subseteq \wpmun$ such that \fr{eigenSig} holds for  any
$\bt$ and $\WW \in \calF_1$, and a set $\calF_2 \subseteq \wpmun$ such that, \eqref{remainder_ineq} hold.  
Let   $\calF = \calF_1 \cap \calF_2$. Lemma \ref{lem:RanTerm} implies that on the event $\mathcal{F}$
\begin{equation*}
\frac{2 \, \left| \langle \hba - \balpha, \boB^T \bbb \rangle \right|}{n} 
\leq \hdel \| \hba - \balpha \|_{\mathrm{block}}.
\end{equation*}

Consider a set $\Xi$ of values of the vector $\bxi$ such that 
\be \label{set_Xi}
 2 \sig  n^{-1} \, \left| \langle \hba - \balpha, \boB^T \bxi \rangle \right| 
\leq \hdel \| \hba - \balpha \|_{block} \quad \mbox{for} \quad \bxi \in \Xi.
\ee 
 Using \eqref{vcm_matr}, for $\bxi \in \Xi$, any $\bt$ and $\WW \in \cal F$, one obtains  
\beqn  \label{eq1}
 \frac{\| \boB (\hba - \ba_0) \|_2^2}{n} \leq   \frac{\| \boB (\balpha - \ba_0) \|_2^2}{n}  
-  2 \hdel  \| \hba  \|_{block}  + 2 \hdel \| \balpha  \|_{block} +  \hdel \| \hba - \balpha \|_{block}.
\eeqn
  Since $\balpha$ is an arbitrary vector, setting $\balpha = \ba$  in \fr{eq1} yields
$$
  \hdel  \| \hba  \|_{block}  -   \hdel \| \ba  \|_{block} \leq 0.5\,  \hdel \| \hba - \ba \|_{block}.
$$ 
Let the set $J_{0}$  contain  the indices of nonzero blocks of $\ba_0$:
\begin{equation}\label{def_set_indices}
J_{0}=\left \{(j,l)\;:\; \Vert\ba_{jl}\Vert_{2}\not=0\right \}.
\end{equation}
Then, the last inequality implies 
\begin{equation}  \label{cone_ineq}
\sum_{(i,j)\in J_0^{C}}\Vert (\ba-\hat \ba)_{ij}\Vert\leq 3\sum_{(i,j)\in J_0}\Vert (\ba-\hat \ba)_{ij}\Vert
\end{equation}

From Lemma \ref{lem:Sigma} it follows that 
 \begin{equation*}
\lamin(\hat \Sigma)\geq (1-h)\phi_{\min}\, \omins \quad\text{and} \quad \lamax(\hat \Sigma) \leq (1+h)\phi_{\max}\, \omaxs 
\end{equation*}
where $\omins$ and $\omaxs$ are defined in \fr{ominmax}. 
For $1 \leq j \leq p$ consider sets 
\beqns 
& G_{00}   = \{j: 1 \leq j \leq p,\ \alpha_{j0} =0 \}, \quad 
& G_{01}   = \{j: 1 \leq j \leq p,\ \alpha_{j0} \neq 0 \}, \\
& G_{j0}   = \{ l:   1 \leq l \leq M,\ \| \balpha_{jl} \|_2 = 0 \},  \quad   
& G_{j1}   = \{l:   1 \leq l \leq M,\ \| \balpha_{jl} \|_2 \neq 0 \}
\eeqns
  We choose $\alpha_{j0} = a_{j0}$ if $a_{j0} \neq 0$ and $\alpha_{j0} = 0$ otherwise.  Let sets $G_{j1}$ be so that $l \in G_{j1}$ iff 
$$\| \ba_{jl}\|_2^2 > \eps = \dfrac{ 8^{2}  \left (\sig^2  C_{\om} K^2\mu+1\right )\log p }{n \lamin (\hbSig)}.$$
We set $  \balpha_{jl} = \ba_{jl}$ if $j \in \calJ$ and $l  \in G_{j1}$ and  
$\balpha_{jl} = \bzero$ otherwise where  $\calJ$ is the set of indices corresponding to non-constant functions   $f_j$.

With $\del = 2\,\hdel$,  inequality \fr{cone_ineq} and Lemma \ref{lemma_bickel_tsybakov} guarantee that  
 \begin{equation}\label{new_1}
 \frac{\| \boB (\hba - \ba) \|_2^2}{n}\geq C_B \lamin(\hat \Sigma)\Vert \hat \ba-\ba\Vert^{2}_{2}.
 \end{equation}
 On the other hand, Lemma \ref{lem:Sigma} and the definition of $\balpha$ imply that 
\begin{equation}\label{new_2}
\frac{\| \boB (\balpha - \ba) \|_2^2}{n}  \leq \lamax(\hat \Sigma) \Vert \balpha-\ba\Vert^{2}_{2}.
\end{equation}  
 Then, using \eqref{new_1} and \eqref{new_2}, we rewrite inequality \fr{eq1} as
\beqn  \label{eq2}
 C_B \lamin(\hat{\Sigma})\,\Vert \hat \ba-\ba_0\Vert^{2}_{2}
   & \leq &  \lamax(\hat{\Sigma})\|  ( \balpha - \ba_0) \|_2^2\\
& + & 4 \hdel \sum_{j \in G_{01}} |\ha_{j0} - a_{j0}| + 4  \hdel \sum_{j \in \calJ} \ \sum_{l \in G_{j1}} \| \hba_{jl} - \ba_{jl} \|_2. \nonumber
\eeqn
Using inequality $2 x_1 x_2 \leq x_1^2 + x_2^2$ for any $x_1$, $x_2$, we derive 
\beqns 
4 \hdel \sum_{j \in \calJ} \ \sum_{l \in G_{j1}}  \| \hba_{jl} - \ba_{jl} \|_2 \leq 
\frac{C_B \lamin (\hbSig)}{2} \sum_{j \in \calJ}  \sum_{l \in G_{j1}}  \| \hba_{jl} - \ba_{jl} \|_2^2
+\sum_{j \in \calJ} \frac{8 \hdel^2 \card(G_{j1})}{C_B \lamin(\hbSig)} ,
\eeqns
and similar inequality applies to the first sum in \fr{eq2}. 
By subtracting  $0.5\, C_B\,\lamin (\hbSig) \,     \| \hba - \ba\|_2^2$ 
from both sides of \fr{eq2} and plugging in the values of $\hdel$ and  $\balpha$,
  derive:
\beqn  \label{eq3}
   \|   (\hba - \ba) \|_2^2 & \leq &  \frac{C_B \, \lamax (\hbSig)}{\lamin(\hbSig)} \ \lkv  \sum_{j=1}^p\  \sum_{l \in G_{j0}}\ \| \ba_{jl} \|_2^2 
+   \frac{ 8^{2} \left (\sig^2  C_{\om} K^2\mu+1\right )    (s_0+s)\log p}{n \lamin (\hbSig)}\right. \\& + &\left .
\sum_{j \in \calJ} \ \frac{ 8^{2} \left (\sig^2  C_{\om}\, K^2\mu+1\right )   \card(G_{j1})\log p}{ n  \lamin (\hbSig)}
\rkv. \nonumber
\eeqn

Observe that for any $1 \leq j \leq p$ one has
\beqns 
 \sum_{l \in G_{j0}}\  \| \ba_{jl} \|_2^2  + 
 \frac{ 8^{2}  \left (\sig^2  C_{\om}\, K^2\mu+1\right ) \card(G_1j)\log p}{ n\, \lamin (\hbSig)} 
\leq    \sum_{l=1}^M \min \lkr  \| \ba_{jl} \|_2^2, \frac{ 8^{2}  \left (\sig^2  C_{\om}\, K^2\mu+1\right )\log p}{n\, \lamin (\hbSig)} \rkr.
\eeqns
Application of inequality \fr{block-error} with $\eps =\dfrac{8^{2}  \left (\sig^2  C_{\om}\, K^2\mu+1\right )\log p }{n \lamin (\hbSig)\log n}$ 
(see Lemma \ref{lem:fun_er}) yields 
\beqn \label{cond-error}
\hskip -1 cm \|    \hba - \ba  \|_2^2 & \leq &    \frac{C_B\, \lamax (\hbSig)}{\lamin(\hbSig)}   
\lkv \frac{  \left (C_{\om}\, \sig^2  K^2\mu+1\right )  (s_0+s)\log p }{n \lamin (\hbSig)} 
 \right. \\
& + & \left. \sum_{j \in \calJ}\   C_a^{2/(2r_j+1)}\   \lkr \frac {C_{\om}\,  \sig^2  K^2\mu+1 }
{n \lamin (\hbSig)}  \rkr^{\frac{2r_j}{2r_j +1}}
(\log n)^\frac{(2-\nu_j)_+-2\nu_j\,r_j}{\nu_j(2r_j +1)}\left (\log p\right )^{\frac{2r_j}{2r_j +1}} \rkv \nonumber
\eeqn

Denote $r^{*}_j = r_j\wedge r'_j$ and $r_{\min}^{*} = \min_j r^{*}_j$.
Choose  $L+1\geq n^{1/2}$ and note that $r_{\min}^{*} \geq 2$.   Using \eqref{tail_bounds}, we obtain 
$\| \hbf - \bof \|^2_2   \leq \|    \hba - \ba  \|_2^2 + C^{2}_{a} s (L+1)^{-2r_{\min}^{*}} $,  which implies 
$$
\| \hbf - \bof \|^2_2   \leq \|    \hba - \ba  \|_2^2 + \dfrac{C^{2}_a\,s}{n^{2}}. 
$$
 Now, \eqref{cond-error} together with Lemmas 1,2 and 4 imply
$$\PP\left (\| \hbf - \bof \|^2_2   \leq \Delta  (s_0, s, n, \br )\right )\geq 1-8\,p^{-\mu}.$$ 

\medskip

{\bf Proof of Theorem  \ref{th:upper_mean}. }
 Let sets $\tcalF$ be such that \fr{cond-error} holds.
Then,  $P(\tcalF)  \geq 1 - 8 p^{-\mu}$ and   \fr{cond-error} yields
\beqns 
\EE \| \hba - \ba \|_2^2 & \leq &  \EE \lkv \| \hba - \ba \|_2^2\, \II(\tcalF) \rkv + 
\EE \lkv \| \hba - \ba \|_2^2\, \II(\tcalF^C) \rkv \\
& \leq & C\, 
\frac{\omaxs}{\omins} \lkv 
 \frac{\sig^2  s_0 \log n}{n \omins} +  \sum_{j \in \calJ}\   C_a^{2/(2r_j+1)}\   \lkr \frac { \sig^2}
 {n \omins }  \rkr^{\frac{2r_j}{2r_j +1}}
 (\log n)^\frac{(2-\nu_j)_+-2\nu_j\,r_j}{\nu_j(2r_j +1)}\left (\log p\right )^{\frac{2r_j}{2r_j +1}} \rkv \\&+ & 16\, C_f^2 \, p^{1-\mu}  \\
& \leq & C\,  \Delta_{upper}  (s_0, s, n, r ),
\eeqns
due to \fr{tau_bound}.


\subsection{Proof of Lemma \ref{lem:Sigma}.}\label{RIC}
      
In order to simplify the notations, we set $\omax(\aleph)=\omax$ and  $\Umu = \Umu(\aleph)$. 
Let $\wpmu$ be the set of points described in condition  \fr{conA3pr} 
of  Assumption {\bf(A5)}.
Denote  the direct product of $n$ sets $\wpmu$ by $\wpmun$.
Then, 
$$
 \PP(\wpmun) \geq 1 - 2\,n\,p^{-2\mu}.
$$
Consider random matrices 
$$
\bZ_i =  (\bSig_i)_{\Lambda}  - \bSig_{\Lambda} = (\bOmega_i)_{\Lambda} \otimes \bPhi_i - \bOmega_{\Lambda} \otimes \bPhi, \quad
\bzeta_i = (\bSig_i)_{\Lambda} \II(\wpmu) - \EE ((\bSig_i)_{\Lambda} \II(\wpmu)).
$$ 
Then, $\bzeta_i$ are   i.i.d. with $\EE \bzeta_i=0$.
We apply  the matrix version of Bernstein's inequality, given in Tropp \cite{tropp-user}:

\begin{proposition} \label{th:tropp} {\bf (Theorem 1.6,  Tropp  (2011))}  
   Let $\bzeta_1,\dots ,\bzeta_n$ be independent random matrices in $\RR^{m_1\times m_2}$
such that $\EE(\bzeta_i)=0$. Define
  \begin{equation*}
  \sigma_{\bzeta} = \max \lfi \left \Vert \dfrac{1}{n}\Sum \EE \left(\bzeta_i \bzeta^{T}_i \right)\right \Vert^{1/2}, 
\left \Vert \dfrac{1}{n}\Sum \EE\left (\bzeta_i^{^{T}}\bzeta_i\right )\right \Vert^{1/2}\rfi.
  \end{equation*}
and suppose that $\Vert \bzeta_i\Vert\leq T$ for some $T>0$. Then,   for all $t>0$, with probability at least $1-e^{-t}$ one has
\be \label{Bernstein} 
\left\Vert \dfrac{1}{n}\Sum \bzeta_i\right\Vert\leq 2\max \left 
\{\sigma_{\bzeta} \sqrt{\dfrac{t+\log(d)}{n}}, T \,\dfrac{t+\log(d)}{n}\right \},
\ee
where $d=m_1+m_2$.
\end{proposition}

 In order to find $\sigma_{\bzeta}$, note that 
\beqns
\left\| \frac{1}{n} \Sum \EE (\bzeta_i \bzeta_i^T) \right \| &  =  & \|  \EE (\bzeta_1 \bzeta_1^T) \|  \leq
\| \EE ((\bSig_1)_{\Lambda} (\bSig_1)_{\Lambda}^T \II(\wpmu) \| +\| \EE ((\bSig_1)_{\Lambda}  \II(\wpmu) \| \| \EE ((\bSig_1)_{\Lambda}^{T}  \II(\wpmu) \| \\
& = & \| \EE\lkv ((\bOmega_1)_{\Lambda} \otimes \bPhi_1)((\bOmega_1)_{\Lambda} \otimes \bPhi_1) \II(\wpmu )\rkv \| 
+ \| \EE\lkv ((\bOmega_1)_{\Lambda} \otimes \bPhi_1)\II(\wpmu )\rkv \| ^{2} \\
&=&
 \| \EE \lkv ((\bOmega_1)_{\Lambda}( \bOmega_1)_{\Lambda}) \otimes (\bPhi_1 \bPhi_1)  \II(\wpmu )\rkv \|
+ \| \EE ((\bOmega_1)_{\Lambda}  \II(\wpmu ))\| ^{2}\| \EE (\bPhi_1)\| ^{2}  \\
& \leq & \| \EE (\bPhi_1 \bPhi_1) \| 
\| \EE \lkv ((\bOmega_1)_{\Lambda} (\bOmega_1)_{\Lambda})   \II(\wpmu) \rkv \|
+ \| \EE ((\bOmega_1)_{\Lambda}  \II(\wpmu ))\| ^{2}\| \EE (\bPhi_1)\| ^{2},   
\eeqns
and, similarly, 
\beqns 
\left\| \frac{1}{n} \Sum \EE (\bzeta_i^T \bzeta_i ) \right \|  & \leq & 
\| \EE (\bPhi_1 \bPhi_1) \|\  \| \EE \lkv ((\bOmega_1)_{\Lambda}( \bOmega_1)_{\Lambda})   \II(\wpmu) \rkv \|
+ \| \EE ((\bOmega_1)_{\Lambda}  \II(\wpmu ))\| ^{2}\| \EE (\bPhi_1)\| ^{2}.
\eeqns
Here, 
$$
\| \EE (\bPhi_1 \bPhi_1) \|  \leq \| C_{\phi}^2 (L+1) \bPhi \| = C_{\phi}^2 (L+1) \phimax
$$
and 
\beqns
 \| \EE \lkv ((\bOmega_1)_{\Lambda}( \bOmega_1)_{\Lambda})   \II(\wpmu) \rkv \|   
&   =  &    \| \EE \lkv \bW_{\Lambda} \bW_{\Lambda}^T \bW_{\Lambda} \bW_{\Lambda}^T \II(\wpmu) \rkv \| \\
& \leq & U_{\mu}^2 \| \EE(\bW_{\Lambda} \bW_{\Lambda}^T) \| = U_{\mu}^2 \| \bOmega_{\Lambda} \| =  U_{\mu}^2 \, \omax,
\eeqns
so that 
\be \label{sig_bzeta}
\sigma_{\bzeta}^2 \leq   2\, C_{\phi}^2\, U_{\mu}^2\, (L+1) \phimax \omax.
\ee
Now, observe that, since matrix $\EE((\bSig_i )_{\Lambda} \II(\wpmu))$ is non-negative definite and matrices 
$\bPhi_i$ and $(\bOmega_i)_{\Lambda}$ have rank  one for any $i$, one has 
\be \label{Tvalue}
T   =   \sup \| \bzeta_1 \|    \leq  2  \sup  \| (\bSig_1)_{\Lambda} \II(\wpmu)) \| 
= 2 \sup   \| (\bOmega_1)_{\Lambda} \II(\wpmu) \| \|\bPhi_1 \| 
\leq 2\,C_{\phi}^2 U_{\mu}^2 (L+1).
\ee
Apply  Bernstein inequality \fr{Bernstein} with $\sigma_{\bzeta}^2$ and $T$ given by formulae \fr{sig_bzeta}
and \fr{Tvalue}, respectively. Then, we  obtain for any $t >0$, 
with probability at least $1-e^{-t}$ 
\be \label{Bern_bzeta}
\left\Vert \dfrac{1}{n}\Sum \bzeta_i \right\Vert \leq  4\, \max \left 
\{ \frac{ C_{\phi}  U_{\mu} \sqrt{ (L+1) (t+\log(L\aleph)\, \phimax \omax )} }{\sqrt{n}},   
\,\dfrac{C_{\phi}^2 U_{\mu}^2 (L+1) (t+\log(L\aleph))}{n}\right \}.
\ee
In order to  apply   inequality \fr{Bern_bzeta} to $\bZ_i$, observe that
$ \bZ_i - \bzeta_i = (\bSig_i)_{\Lambda} \II(\wpmu^c) - \EE((\bSig_i)_{\Lambda} \II(\wpmu^c)$ and 
\beqn \nonumber
\| \EE ((\bSig_i)_{\Lambda} \II(\wpmu^c)) \|^2  & = & \| \EE \left ((\bOmega_i)_{\Lambda} \II(\wpmu^c)\right ) 
\otimes \EE (\bPhi_i)\|^2\leq \| \EE \left ((\bOmega_i)_{\Lambda} \II(\wpmu^c)\right )\|^2 \| \EE (\bPhi_i)\|^2 
  \\
& \leq &   \left (\EE \| (\bOmega_i)_{\Lambda} \II(\wpmu^c)\|_{2}\right )^2\, \EE \| \bPhi_i \|^{2}_2 
\leq 2\, C_{\phi}^4\, V \aleph\, p^{-2\mu}  (L+1)^{2},
\label{expec_bound}
\eeqn 
due to the fact that $\EE (\bW^{(j)})^4 \leq V$  for any $j=1, \cdots,p$.

Now, we use the union bound over all $\Lambda$ such that $\vert\Lambda\vert\leq \aleph$, 
the inequality $\binom{p}{\aleph}\leq \left (\frac{ep}{\aleph}\right )^{\aleph}$ 
and choose $t=2\,\mu\,\aleph\,\log(\frac{ep}{\aleph})$.
Combining \fr{Bern_bzeta} and \fr{expec_bound}, for all $\Lambda$ such that $\vert\Lambda\vert\leq \aleph$, we derive  
\begin{align*}
&  \inf_{ \Lam:\, | \Lam |  \leq \aleph }\  \PP \lkr  \lfi \left\| \frac{1}{n} \Sum \bZ_i \right\| < z \rfi \cap \lfi \WW \in \wpmun \rfi \rkr \\
& \geq 
 \inf_{ \Lam:\, | \Lam |  \leq \aleph }\  \PP \lkr  \lfi \left\| \frac{1}{n} \Sum \bzeta_i \right\| < z - \| \EE ((\bSig_1)_{\Lambda} \II(\wpmu^c)) \| \rfi 
\cap \lfi \WW \in \wpmun \rfi \rkr \\
& \geq 
1 -  \sup_{ \Lam:\, | \Lam |  \leq \aleph }\  \PP \lkr \lfi \left\| \frac{1}{n} \Sum \bzeta_i \right\| 
\geq z -  C^{2}_{\phi}\,(L+1)\,  p^{-\mu}\sqrt{ 2\,V \, \aleph }  \rfi \cap \lfi \WW \in \wpmun \rfi \rkr 
\geq 1 -  (ep)^{-\mu}- 2\,n\,p^{-2\mu }
\end{align*}
for any $z$ such that
\beqn 
z & \geq &   8\, \max \left 
\{ \frac{ C_{\phi}  U_{\mu} \sqrt{\mu\,\aleph\, (L+1) \log (p+L)\, \phimax \omax   }}{\sqrt{n}},   
\  \dfrac{\mu\,\aleph\,C_{\phi}^2 U_{\mu}^2 (L+1)  \log (p+L)}{n}\right \} \hspace{4mm}
\nonumber\\
& + &   C^{2}_{\phi}\, (L+1)\, \,  p^{-\mu}\sqrt{2 V  \, \aleph}.
 \label{t_cond}
\eeqn
Note that, under condition \fr{mu_cond}, one has
$$
 C^{2}_{\phi}\,(L+1)\, \,  p^{-\mu}\sqrt{ 2\,V \, \aleph }  \leq8\,\mu\aleph C_{\phi}^2 U_{\mu}^2 (L+1)\, \log(p+L)\   n^{-1}.
$$
It is easy to check that,  whenever   $n \geq N (\lam)$ where $N (\lam)$ is defined in \fr{No},
condition \fr{t_cond} is satisfied with  
$$ 
z =  \dfrac{8 C_{\phi}  U_{\mu} \sqrt{\mu\,\aleph\, (L+1) \phimax \omax  \log (p+L) }}{\sqrt{n}}  + 
\dfrac{9 \, \mu\,\aleph\,C_{\phi}^2 U_{\mu}^2 (L+1)  \,\log(p+L)}{n} \leq h \,\omin\, \phimin,
$$
which, together with condition $p^{\mu}\geq 2n$, implies that 
\be \label{ld_ineq}
\inf_{ \Lam:\, | \Lam |  \leq \lam }\  \PP\lkr \|\hbSig_{\Lambda} - \bSig_{\Lambda}\| \leq h \omin \phimin \rkr \geq 1 -  p^{ - \mu}-p^{-\mu}.
\ee
In order to complete  the proof, observe that  $\lamin(\hbSig_{\Lambda}) \geq \lamin(\bSig_{\Lambda}) - \|\hbSig_{\Lambda} - \bSig_{\Lambda}\|$ and
 $\lamax(\hbSig_{\Lambda}) \leq \lamax(\bSig_{\Lambda}) + \|\hbSig_{\Lambda} - \bSig_{\Lambda}\|$. 
 \\


\subsection{Proofs of supplementary lemmas}
\label{sec:proofs-supplementary}


\begin{lemma}\label{lem:RanTerm}
Let vector $\balpha \in \RR^{p(L+1)}$ be partitioned into subgroups the way 
it has been done for  vector $\ba = \Vect(\bA)$. Let $L+1\geq n^{1/2}$ and $n\geq \mathbf N$
where $\mathbf N$ is defined in \fr{Nmax}. 
 Suppose that Assumptions {\bf (A5)}, {\bf (A4)} are valid and  $\mu$ in \fr{conA3pr} is large enough, 
so that condition \eqref{mu_cond_2} holds.
 Then, \begin{itemize}
\item[(i)]
\be \label{LargeDevXi}
\PP \lkr \frac{2 \sig  \, \left| \langle \balpha, \boB^T \bxi \rangle \right|}{n} 
\leq \hdel \| \balpha \|_{\mathrm{block}}  \rkr \geq 1 - 5p^{-\mu};
\ee
\item[(ii)]\be \label{LargeDevB}
\PP \lkr \frac{2 \, \left| \langle \balpha, \boB^T \bbb \rangle \right|}{n} 
\leq \hdel \| \balpha \|_{\mathrm{block}}  \rkr \geq 1 -   4p^{-\mu},
\ee
\end{itemize}
where $\hdel$ is defined in \fr{deldef}.
 \end{lemma}


{\bf Proof.}  In order to prove (i), note that 
\beqn \label{inequa_1_random}
n^{-1/2}\, \left| \langle \balpha, \boB^T \bxi \rangle \right| & \leq &
n^{-1/2}\,  \max_{1 \leq j \leq p} |(\boB^T \bxi)^{(j,0)}| \  \sum_{j=1}^p |\alpha_{j0}| \\
& + &
n^{-1/2}\, \max_{\stackrel{1 \leq j \leq p}{1 \leq l \leq M}} \sqrt{\sum_{k \in K_l} \lkv (\boB^T \bxi)^{(j,k)} \rkv^2}  \   
\sum_{j=1}^p \sum_{l=1}^M  \sqrt{\sum_{k\in K_l}\alpha^{2}_{jk}}. 
\nonumber
\eeqn

  Fix vectors $\bW_1, \cdots, \bW_n$ and $\bt$.  Using Hoeffding-type inequality for sub-gaussian random variables 
(see, e.g., Proposition 5.2 in \cite{vershynin}) we obtain that, for any $\mu \geq 1$,
\be \label{prob1_random}
\PP \lkr n^{-1/2}\,  \max_{1 \leq j \leq p} |(\boB^T \bxi)^{(j,0)}| \leq 
\sqrt{C\,K^2 \mu \, \aleph_{\max}(\hbSigma_{\Lambda})\,\log(p) }\ \Bigg |   \bt, \WW   \rkr
\geq 1 - e\, p^{1-2\mu}
\ee
where $C>0$ is an absolute constant and $\vert \Lambda\vert\leq 1$.

For the second maximum in \fr{inequa_1_random},   we use a corollary of the Hanson-Wright inequality for sub-gaussian random vectors  
(see Theorem 2.1 in \cite{rudelson_vershynin}), which states that, for any matrix $\bQ\in \mathbb R^{n\times d}$ and any $\mu\geq 1$,  one has
$$
\PP \lkr \| \bQ \bxi \|_{2}  \leq \Vert \bQ\Vert_{2} + C\,K \lamax(\bQ) \sqrt{\mu\log(p)}\rkr \geq 1 - 2p^{-3\mu}. 
$$
We apply this inequality with  matrix $\bQ =  n^{-1/2} \,\boB_{j,l}$,    where $\boB_{j,l}$ is a  
$n \times d$ sub-matrix of matrix $\boB$. Note  that $\lamax(\bQ) \leq  \sqrt{\lamax(\hbSigma_{\Lambda})}$ where 
$\vert \Lambda\vert\leq 1$.  Using $\Vert \bQ\Vert^2_{2} \leq d \lamax(\hbSigma_{\Lambda})$ and $d \approx \log n$, obtain
for any $\mu \geq 1$  
\be \label{prob2}
\PP \lkr n^{-1/2}\, \max_{\stackrel{1 \leq j \leq p}{1 \leq l \leq M}}  {\sum_{k \in K_l} \lkv (\boB^T \bxi)^{(j,k)} \rkv 
 \leq 
 (1+ CK\sqrt{\mu})\sqrt{\lamax (\hbSigma_{\Lambda})\,\log (p)}  } \Bigg |   \bt, \WW   \rkr
\geq 1 -   \frac{2\,p L}{   p^{3\mu}\, \log (n)  }
\ee
 where we used $M=L/\log (n)$.  
 
 Now, consider a set  $\calF_1 \subseteq \wpmun$ such that  \fr{eigenSig} holds for  any
 $\bt$ and any $\WW \in \calF_1$.  By Lemma \ref{lem:Sigma}, obtain $\mathbb P(\mathcal{F}_{1})\geq 1-2\,p^{-\mu}$. Then,   
 \beqn  
 n^{-1}\, \left| \langle \balpha, \boB^T \bxi \rangle \right| & \leq & CK\sqrt{\dfrac{\mu (1+h)\phimax\omax(1)\log p}{n}}
 \nonumber
 \eeqn
 with probability at least $1-5\,p^{-\mu}$, due to condition \eqref{mu_cond_2}.

Inequality   (ii) is proved  in a similar manner. Indeed, one has  
\beqn \label{inequa_1}
n^{-1/2}\, \left| \langle \balpha, \boB^T \bbb \rangle \right| & \leq &
n^{-1/2}\,  \max_{1 \leq j \leq p} |(\boB^T \bbb)^{(j,0)}| \  \sum_{j=1}^p |\alpha_{j0}| \\
& + &
n^{-1/2}\, \max_{\stackrel{1 \leq j \leq p}{1 \leq l \leq M}} \sqrt{\sum_{k \in K_l} \lkv (\boB^T \bbb)^{(j,k)} \rkv^2}  \   
\sum_{j=1}^p \sum_{l=1}^M  \sqrt{\sum_{k\in K_l}\alpha^{2}_{jk}}. 
\nonumber
\eeqn
Consider a set  $\calF_1 \subseteq \wpmun$ such that  \fr{eigenSig} holds for  
any  $\bt$ and any $\WW \in \calF_1$.   By Lemma \ref{lem:Sigma},  obtain that 
$\mathbb P(\mathcal{F}_{1})\geq 1-2p^{-\mu}$. Then, on the event $\mathcal{F}_{1}$, the following inequalities are valid
\begin{equation}\label{RIC_1}
 \begin{split}
 n^{-1/2}\max_{1 \leq j \leq p} |(\boB^T \bbb)^{(j,0)}|& \leq \sqrt{(1+h)\phi_{\max}\, \omax(1)}\ \Vert \bbb\Vert_{2},  \\
 n^{-1/2}\max_{\stackrel{1 \leq j \leq p}{1 \leq l \leq M}} \sqrt{\sum_{k \in K_l} \lkv (\boB^T \bbb)^{(j,k)} \rkv^2} 
& \leq \sqrt{(1+h)\phi_{\max}\, \omax(1)}\ \Vert \bbb\Vert_{2}.
 \end{split}
\end{equation}

  Now consider a set $\calF_2 \subseteq \wpmun$  such that \eqref{remainder_ineq} holds. Let   $\calF = \calF_1 \cap \calF_2$.
  Lemmas \ref{lem:Sigma} and \ref{lem:remainder} imply that on the event $\mathcal{F}$
 \begin{equation*}
  \begin{split}
  n^{-1}\max_{1 \leq j \leq p} |(\boB^T \bbb)^{(j,0)}|& \leq   C_a   \, (L+1)^{-2}
\sqrt{3\,g_2\,s\,(1+h)\phimax\, \omax(s)\, \omax(1)} \quad \text{and}\\
 n^{-1} \max_{\stackrel{1 \leq j \leq p}{1 \leq l \leq M}} \sqrt{\sum_{k \in K_l} \lkv (\boB^T \bbb)^{(j,k)} \rkv^2}
& \leq     C_a   \, (L+1)^{-2}\sqrt{3\,g_2\,s\,(1+h)\phimax\, \omax(s)\, \omax(1)}. 
   \end{split}
  \end{equation*}
  Now,   the statement of the Lemma is valid due to assumption $(L+1)^2 \geq n$.

\medskip


\begin{lemma}\label{lemma_bickel_tsybakov} 
Let $\mathcal{S}\in\{1,\dots,p\}$ be a subset of indices with $\vert \mathcal{S}\vert\leq (s+s_0)$ and 
 $J =S\times \{0,\dots,L\}\subset \{1,\dots,p\}\times\{0,\dots, L\}$. 
 We define $\mathcal{A}$ to be the following set of vectors
\begin{equation*}
\mathcal{A}=\left \{\bv\in \mathbb R^{(L+1)p}\;:\;\sum_{(i,j)\in J^{C}}\Vert \bv_{ij}\Vert\leq 3\sum_{(i,j)\in J}\Vert \bv_{ij}\Vert\right \}.
\end{equation*}

Suppose that the assumptions of Lemma \ref{lem:Sigma} and Assumption {\bf (A6)} hold. Then, there exists a numerical constant $C_B$ such that, 
\begin{equation}\label{cond_bickel_tsybakov}
\underset{\mathcal{S}}{\min}\;\underset{\bv\in \mathcal{A}}{\min} \dfrac{\Vert \boB\bv\Vert}{\Vert \bv\Vert}\geq C_B (1-h)\phi_{\min}\, \omins
\end{equation}
 where $\omins$ is defined in \fr{ominmax}.

\end{lemma}

\noindent
{\bf Proof.} The proof follows the lines of the proof of  Lemma 4.1 in \cite{bickel_ritov_tsybakov}.

\medskip
 

 \begin{lemma} \label{lem:fun_er}
Let function $f(t) = \sum_{k=1}^\infty a_k \phi_k(t)$ satisfy condition {\bf (A4)}, i.e. 
\be \label{confun}
\sum_{k=0}^\infty |a_{k}|^\nu (k+1)^{\nu r '} \leq C_a^\nu, \quad r ' = r  +1/2-1/\nu,
\ee
 for some $C_a>0$, $1 \leq \nu < \infty$ and $r  > \min(1/2, 1/\nu)$.
Let $\ba_l$ be blocks of coefficients $a_{k}$ of length $d$, so that $\ba_l = (a_{(l-1)d+1}, \cdots, a_{ld})$.
Then, for any $\eps >0$ one has
\be \label{block-error}
\sum_{l=0}^\infty \min(\| \ba_l \|_2^2, \eps  d ) \leq   C_a^{\frac{2}{2r+1}} 
\eps^{\frac{2r}{2r+1}} d^{\frac{(2-\nu)_+}{\nu(2r+1)}}.
\ee
Moreover, if $r' \geq 2$ and basis $\{ \phi_k\}$ satisfies assumption  {\bf (A1)},
then 
\be \label{uniform_bound}
\| f - f_J \|_\infty \leq C_a  C_\phi   (J+1)^{-(r^*-1/2)}\quad \mbox{with} 
\quad f_J (t) = \sum_{k=1}^J a_k \phi_k(t).
\ee
 Here, $r^* = \min(r,r')$, $(x)_+ = x$ if $x >0$ and zero otherwise.
\end{lemma}

{\bf Proof.} First, let us show that for $J \geq 1$
\be \label{tail_bounds}
 \sum_{k=J+1}^\infty a_k^2 \leq    C_a^2 (J+1)^{-2r^*}.
\ee
 Indeed, if $\nu \geq 2$, the Cauchy inequality yields
\beqns
\sum_{k=J+1}^\infty  a_k^2 & \leq & \lkr \sum_{k=J+1}^\infty    |a_k|^\nu k^{r'\nu} \rkr^{2/\nu} 
\lkr \sum_{k=J+1}^\infty k^{- \frac{2r'\nu}{\nu-2}} \rkr^{1 - 2/\nu} \leq
C_a^2 \lkr \frac{\nu-2}{2r\nu} \rkr^{\nu/2-1} (J+1)^{-2r}.
\eeqns
Since $(\nu -2)/(2 r \nu) < 1$ for $\nu \geq 2$, inequality  \fr{tail_bounds} holds. 
If $1 \leq \nu <2$, then 
\beqns 
\sum_{k=J+1}^\infty a_k^2  &\leq& \lkr \max_{k \geq J+1} |a_k|  \rkr^{2-\nu} (J+1)^{-r' \nu} 
\lkr \sum_{k=J+1}^\infty    |a_k|^\nu k^{r'\nu} \rkr \leq
C_a^2    (J+1)^{-2r'},
\eeqns
so that \fr{tail_bounds} is valid.

Now, using \fr{tail_bounds}, we prove \fr{block-error}. Again, we consider cases $\nu \geq 2$ 
and $1 \leq \nu <2$, separately.
If $\nu \geq 2$, then partitioning the sum into the portion for  $l \leq J$  and $l> J$ (which corresponds to $k > Jd$), we derive
\beqns
\sum_{l=1}^\infty \min(\| \ba_l \|_2^2, \eps  d ) & \leq & J \eps d + C_a^2  (J d)^{-2r}.
\eeqns
Minimizing the last expression with respect to $J$, we obtain   \fr{block-error} without the log-factor.
If $1 \leq \nu <2$, then 
\beqns
\sum_{l=1}^\infty \min(\| \ba_l \|_2^2, \eps  d ) & \leq & J \eps d + (d \eps)^{1 - \nu/2} 
\sum_{l= J+1}^\infty  \|\ba_l\|_2^\nu. 
\eeqns
Since for  $1 \leq \nu <2$  
\beqns
\|\ba_l\|_2^2 = \sum_{k=(l-1)d+1}^{ld} a_k^2 & \leq & \lkr \sum_{k=(l-1)d+1}^{ld} |a_k|^\nu \rkr^{2/\nu},
\eeqns
one has 
\beqns
\sum_{l= J+1}^\infty  \|\ba_l\|_2^\nu \leq \sum_{k = Jd +1}^\infty |a_k|^\nu & \leq & C_a^\nu (Jd+1)^{-r'\nu}
\eeqns
and 
\beqns
\sum_{l=1}^\infty \min(\| \ba_l \|_2^2, \eps  d ) & \leq & 
J \eps d +    C_a^\nu (d \eps)^{1 - \nu/2}  (Jd)^{-r'\nu}.
\eeqns
Minimization of the last expression with respect to $J$ yields \fr{block-error}.

In order to prove \fr{uniform_bound}, observe that for any $J > 1$  and $r^* > 3/2$, one has
\beqns
\| f - f_J \|_\infty & \leq &  \underset{t\in [0,1]}{\sup} \sum_{l=1}^\infty \sqrt{ \sum_{k=Jl}^{J(l+1)-1} a_k^2}
\sqrt { \sum_{k=Jl}^{J(l+1)-1} \phi_k^2(t)} 
\leq   C_a C_\phi  \sum_{l=1}^\infty \sqrt{J l} (J l+1)^{-r^*} \\
&\hskip 1 cm \leq &  C_a C_\phi J^{-(r^*-1/2)}
\eeqns
which completes the proof.
\\

\medskip


\begin{lemma} \label{lem:remainder}
Let  $r^{*}_j = r_j\wedge r'_j$, $r^{*} = \min_j r^{*}_j  \geq 2$   in assumption {\bf(A3)}.
 Let $\mu$ in \fr{conA3pr} be large enough, so that 
   \be \label{mu_cond_1}
   p^{\mu} \geq 2\,n
   \ee
and $n$ be such that   
\begin{equation}\label{cond_n_remaind}
n\geq \dfrac{U_{\mu}^{2}\,C_{\phi}^{2}(L+1)\mu\log p}{g_2\, \omax (s)} 
\end{equation}
where $U_{\mu} = U_{\mu} (s + s_0)$.  Then,  one has 
\be \label{remainder_ineq}
\PP \lkr  \lfi n^{-1} \| \bbb \|_2^2\leq  3\,g_2\,\omax(s)\,s\,C_a^{2}(L+1)^{-4} \rfi \cap \lfi \WW \in \wpmun\rfi  \rkr 
\geq 1 -2\,p^{-\mu}.
\ee
Here $\bbb$ is the vector with components  $\bbb_i = \bW_i^T \brho (t_i)$, $i = 1, \cdots, n$,
where  $\brho(t) = (\rho_1(t), \cdots, \rho_p(t))^T$ and $\rho_i(t)$ are  defined in \fr{func_expan}. 
\end{lemma}

{\bf Proof.}  In order to estimate $\| \bbb \|_2^2$ we apply  Bernstein inequality  to the centered random variables 
$\beta_i = b_i^2\II(\wpmu) - \EE \left (b_i^2\II(\wpmu)\right )$. We start with establishing an upper bound for  $\EE  b_i^2$:  
\begin{equation}\label{lemma5_1}
\begin{split}
\EE  b_i^2 & = \EE  \left (\,n^{-1}\Vert \bbb\Vert^{2}_{2}  \right )=
\EE\langle \brho(t)\,\bOmega_{1} \brho(t)\rangle\\&\leq  \omax(s) \EE\Vert \brho(t)\Vert^{2}_{2}.
\end{split}
\end{equation}
Using  \eqref{tail_bounds} and the orthonormality of the basis, we derive
\begin{equation}\label{lema5_3}
\EE\Vert \brho(t)\Vert^{2}_{2}
\leq s\,g_2\, C_a^{2}(L+1)^{-2r^{*}}.
\end{equation}
Plugging \eqref{lema5_3} into \eqref{lemma5_1}, we obtain
\beqns \label{remain_mean}
\EE b_i^2 &  \leq g_2\,\omax(s)\,s\,C_a^{2}(L+1)^{-2r^{*}}.
\eeqns
Now, we use the upper bound \eqref{uniform_bound} to establish an upper bound for the variance $\sigma_b^2$ of $b_i^2$: 
 \beqn 
\label{var_rem}
\begin{split}
\sigma_b^2 
& \leq  \EE \left (b_i^4 \II(\wpmu)\right )\leq   U_{\mu}^{4}\,\EE\Vert \brho(t)\Vert^{4}_{2}\leq g_2\,U_{\mu}^{4}\,s^{2}\,C_a^4  C_\phi^{4}(L+1)^{-4r^{*}+2}.
\end{split}\eeqn
It is also easy to see that, for any $i$ and $\WW \in \wpmun$, by  \fr{uniform_bound}, one has, 
\beqn \label{unif_rem}
U_b = \max (b_i^2) \leq \Vert\bW_{J}\Vert^{2}_{2} \max_t \| \rho(t) \|_2^2\ \leq  
\Umu^2\,C_a^2 C_\phi^2   (L+1)^{-(2r^{*} -1)}\ s.
\eeqn

 With $\sigma^{2}_{b}$ and $U_{b}$ given by \eqref{var_rem} and \eqref{unif_rem}, respectively,  and $t=\mu \log(p)$, 
one obtains, that with probability at least $1-p^{-\mu}$     
\be \label{Bern_remain}
\left\vert \dfrac{1}{n}\Sum \beta_i \right\vert \leq 2\max \left 
\{  \sigma_{b}\sqrt{ \dfrac{\mu \log p}{n}},   
\,U_b\dfrac{\mu \log p}{n}\right \}.
\ee
Since $b^{2}_i- \beta_i=b^{2}_{i}\II(\wpmu^c)+\EE\left (b^{2}_{i}\II(\wpmu)\right )$,  we  derive

\begin{align*}
& \PP \lkr  \lfi  n^{-1} \| \bbb \|_{2}^2  < z \rfi \cap \lfi \WW \in \wpmun \rfi \rkr 
\geq 
\PP \lkr  \lfi \left\vert \dfrac{1}{n}\Sum \beta_i \right\vert < z - \EE  \lkr n^{-1} 
\| \bbb\|_2^2 \,\II(\wpmu)\rkr \rfi \cap \lfi \WW \in \wpmun \rfi \rkr \\
& \geq \PP \lkr \lfi \left | \frac{1}{n} \Sum \beta_i \right | 
< z - g_2\,\omax(s)\,s\,C_a^{2}(L+1)^{-2r^{*}}\rfi 
\cap \lfi \WW \in \wpmun \rfi \rkr
\geq 1 -p^{-\mu}-2\,n\,p^{-2\mu} \geq 1 -2p^{-\mu}
\end{align*}
for any $z$ such that
\beqns 
z  \geq   2\max \left 
\{  \sigma_{b}\sqrt{ \dfrac{\mu\log p}{n}},   
\,U_b\dfrac{\mu \log p}{n}\right \}
+ g_2\,\omax(s)\,s\,C_a^{2}(L+1)^{-2r^{*}}.
\eeqns
For $n$, satisfying condition  \eqref{cond_n_remaind}, one can choose 
$$
z=3\,g_2\,\omax(s)\,s\,C_a^{2}(L+1)^{-2r^{*}}
$$
which together with $r^{*}\geq 2$ implies the statement of the Lemma.
\\



\end{document}